\documentclass[twoside, reqno]{amsart}

\usepackage[left=2.5cm, right=2.5cm, top=3cm, bottom=2.5cm]{geometry}

\usepackage[colorlinks=true, pdfstartview=FitV, linkcolor=blue,citecolor=blue, urlcolor=blue]{hyperref}

 \usepackage{colortbl}
\usepackage[dvipsnames]{xcolor}
\colorlet{palegolden}{Gray!10!white}
\colorlet{palesilver}{Purple!10!white}
\definecolor{labelkey}{rgb}{0,0,1}
\definecolor{Red}{rgb}{0.7,0,0.1}
\definecolor{Green}{rgb}{0,0.7,0}
\usepackage{mymacros}

\usepackage{verbatim} \usepackage{marginnote}
\usepackage{amsfonts}
\usepackage{graphicx} \usepackage{setspace} 
\usepackage{enumerate} \usepackage{mathrsfs} \usepackage{empheq}
\usepackage{graphics} \usepackage{epsfig} 
\usepackage{stackrel}
\usepackage{color}
\usepackage{comment}

\newcommand{\vincent}[1]{{#1}}

\begin{document}

\newtheorem*{main}{Theorem}
\newtheorem{claim}{Claim} 

\title{A PDE model for chemotaxis with logarithmic sensitivity and logistic growth}
\author{Padi Fuster Aguilera, Vincent R. Martinez, Kun Zhao}

\date{December 15, 2020}

\maketitle

\begin{abstract}
    In this paper, we study the initial-boundary value problem \vincent{and its asymptotic behavior} for a repulsive chemotaxis model with logarithmic sensitivity and logistic growth. We \vincent{establish global} well-posedness of \vincent{strong} solutions for large initial data with Neumann boundary conditions and, \vincent{moreover, establish the qualitative result that both the population density and chemical concentration asymptotically converge to constant states with the population density specifically converging to its carrying capacity.} \vincent{We additionally prove} that \vincent{the vanishing chemical diffusivity limit holds in this regime. Lastly,} we provide \vincent{numerical confirmation of the rigorous qualitative results, as well as numerical simulations that demonstrate a separation of scales phenomenon}. 
\end{abstract}



\section{Introduction}
Chemotaxis is the movement of an organism in response to a chemical stimulus. This phenomenon occurs in a variety of ways, for example, with sperm swimming towards the egg during fertilization, bacteria being attracted by food, cancer metastasis, etc.
A PDE model was first derived by Patlak \cite{Patlak53}
in the 50's and later by Keller-Segel \cite{K-S} in the 70's; it has since been studied along with several variants. {\vincent{This system can generally be written as}}
\begin{align}\label{eq:keller:segel}
    \begin{split}
    u_t&=-\nabla\cdot(-D\nabla u+\chi(u,c)
    \nabla c)+f(u,c),\\
    c_t&=-\nabla\cdot(-\veps\nabla c)+g(u,c),
    \end{split}
\end{align}
where $u$ \vincent{denotes the population density} of the organism, $c$ the \vincent{chemical concentration}. {\vincent{The chemotactic flux, as represented by $\chi(u,c)\nabla c$, is a characteristic feature of this type of movement. Here, $\chi$ regulates the sensitivity of the population to the chemical density and the prefactor}}. {\vincent{The parameter, $D$, denotes}} the diffusion coefficient of $u$, \vincent{while $\veps$ denotes the chemical diffusivity.}

In this paper, we will study a particular \vincent{modification of \eqref{eq:keller:segel}, namely, one that assumes logarithmic sensitivity \cite{stevens1997aggregation}, non-linear growth for the chemical concentration, and also accounts for the effects of logistic growth. In particular, we will assume that}
    \begin{align}\label{eq:assumptions}
      \chi(u,c)=\chi \frac{u}c,\quad f(u,c)=\kappa_1u\left(1-\frac{u}{\kappa_2}\right),\quad g(u,c)=\mu uc-\sigma c,
    \end{align}
\vincent{where $\kap_1,\kap_2>0$, $\mu,\s,\chi\in\mathbb{R}$; $\chi$ is attractive when $\chi>0$ and repulsive when $\chi<0$, while $|\chi|$ \vincent{measures} the strength of chemotactic sensitivity. Notice that under these assumptions, $\chi(u,c)\nabla c=\chi u\nabla \phi(c)$, where $\phi(c)=\ln c$}. The \vincent{choice of a} logarithmic (singular) sensitivity, $\phi$, follows \vincent{from} the Weber-Fechner law, which states that subjective sensation of $c$ by $u$ can be given as a proportion of the growth over the density itself. Mathematically, it is a major source of analytical difficulties for obtaining qualitative results. \vincent{The non-linear function, $g$, contains the reaction terms of the system and dictates the growth and death rates of $c$.} Our model also \vincent{accounts for the effects} of logistic growth for the organism's density, \vincent{as modeled by $f$;} this \vincent{was} first omitted by Keller and Segel \cite{K-S} because of the time-scale difference between the movement of the organism and its \vincent{growth and death} rates. However, for certain phenomena, such as organisms moving very slowly, e.g., metastatic cells or organisms moving in semi-solid medium, e.g., in agar or mayonnaise, it is relevant to take into an account the growth of the population (cf. \cite{Mur1}, \cite{Mur2}). 

\vincent{We will specifically consider \eqref{eq:keller:segel} under the assumptions \eqref{eq:assumptions} over the bounded domain $\Om=(0,1)$ complemented by Neumann boundary conditions for $(u,c)$. By applying the standard Cole-Hopf transformation for this system, i.e., setting $v=(\ln c)_x$, as in \cite{zeng2019logarithmic}, for instance, we arrive at}
\begin{align}
    u_t+\chi(u v)_x&=D u_{xx}+ \kappa_1u\left(1-\frac{u}{\kappa_2}\right),\quad u_x(0,t)=u_x(1,t)=0\label{eq:main1}\\
    v_t+\left(\mu u -\veps v^2\right)_x&=\veps v_{xx},\quad v(0,t)=v(1,t)=0.\label{eq:main2}
\end{align}
\vincent{For our purposes, we will assume that $\chi\mu>0$; this guarantees that the system is hyperbolic}. In the case where $\chi\mu<0$, Levine-Sleeman \cite{levinesleeman} demonstrate that solutions can blow-up in finite time. \vincent{As we are ultimately interested in the setting where the chemical diffusivity is small relative to the diffusion of the organism, we will also assume that $\veps/D\ll1$ . Upon non-dimensionalizing the system accordingly, we obtain}
\begin{align}
    u_t+(u v)_x&= u_{xx}+ ru(1-u),\quad u_x(0,t)=u_x(1,t)=0\label{eq:main11}\\
    v_t+\left(u -\frac\veps\chi v^2\right)_x&=\frac\veps D v_{xx},\quad v(0,t)=v(1,t)=0.\label{eq:main22}
\end{align}
\vincent{where $r={\kappa_1 D}/({\chi \mu\kappa_2})$. We will lastly assume that $\chi=D=1$;} biologically, the scale of $\veps$ compared with the ones of $\chi$ and $\mu$ are in accordance with this assumption in many phenomena, for example, in tumor angiogenesis (cf. \cite{sleeman2001partial}). 

\vincent{In contrast with studies of the model without logarithmic growth \cite{li2012global}, which proved that the organism density converges to its initial average, we ultimately prove that incorporation of logistic growth can qualitatively change the long-time behavior of the population density by forcing it to} converge to its carrying capacity instead. \vincent{Systematic numerical simulations are carried out to confirm this qualitative behavior, as well as to expose a separation of scales phenomenon between the diffusive and logistic time-scales. A related work by Zeng and Zhao in \cite{zeng2019logarithmic} also studied \eqref{eq:keller:segel} under the assumptions \eqref{eq:assumptions}, namely, with logistic growth, but} on the whole real line. One of the main difficulties \vincent{encountered in studying this model on bounded domains, as dealt with in this paper,} is given \vincent{by} the lack of conservation of mass in the chemical concentration equation. \vincent{Indeed, this feature is relied on in \cite{zeng2019logarithmic}. As a result, deriving estimates analogous to \cite{zeng2019logarithmic} is simply not a straightforward adaptation to bounded domains. Instead, we identify a Lyapunov structure (\cref{sect:Lyapunov}) in our system. By additionally taking advantage of the compatibility conditions at the boundary, we successfully overcome the apparent obstruction of the lack of mass conservation}. 

\vincent{The remainder of the paper is organized as follows: First, we establish the mathematical setting and notation that we use throughout the paper in \cref{sect:background}. We then develop the underlying Lyapunov structure of our system (see \cref{sect:Lyapunov}); this is a crucial feature for obtaining global existence and uniqueness of strong solutions for the initial-boundary value problem on $\Omega=(0,1)$ with Neumann boundary conditions for the transformed system $(u,v)$, where $v=\partial_x \ln c$ (see \cref{thm:Glob:exist}), which is proven in \cref{sect:Glob:exist}. With the Lyapunov structure and global existence theory in hand, we then prove the zero-chemical-diffusivity limit holds in the topology of uniform convergence. We, in fact, show the convergence occurs in a stronger topology, namely, that of the Sobolev space $H^1$ (see \cref{thm:diff}). We conclude the paper with some numerical results that verify the qualitative result and moreover, expose a potential relation between the time scales of the chemical diffusive parameter and the logistic growth that suggests further investigation. In a follow up paper, we will address the case of Dirichlet boundary conditions, as well as treat the issue of boundary layers (singular limit) in the corresponding zero-chemical diffusivity problem.}

\section{Mathematical Background}\label{sect:background}

In this section, we define the functional setting in which we work, as well as the notation and conventions that we use for the rest of the paper. \vincent{From now on, we will suppose $\Omega=(0,1)$.}


The $L^p$-based Sobolev spaces of order $k$, over $\Omega$, where $1\leq p\leq\infty$ and $k>0$ is defined by
\begin{align}\notag
&\vincent{H^{k}(\Omega):=\{h\in L^2(\Omega): D^\alpha h\in L^2(\Omega), \text{for all}\ |\alpha|\leq k\},\quad D:=\frac{d}{dx}.}
\end{align}
\vincent{We also define }
    \begin{align}\notag
        &\vincent{H^{1}_0(\Omega):=\overline{C^\infty_c(\Omega)}^{H^{1}},}
    \end{align}
\vincent{where $C^\infty_c(\Omega)$ denotes the space of smooth, compactly supported test functions over $\Omega$ and the bar denotes the closure with respect to the norm $\lVert\cdotp\rVert_{H^{1}}$.} 
\begin{align}
    \vincent{\|h \|^2}&=\| h \|_{L^2(\Omega)}^2=\intol h(x)^2 dx\notag\\
      \vincent{ \|h\|_p^p}&=\| h \|_{L^p(\Omega)}^p=\intol \vincent{|h(x)|^p dx}\notag\\
    \|h\|_\infty&=\|h\|_{L^\infty(\Omega)}= \esssup_{x\in\Omega}|h(x)|\notag\\
      \vincent{  \Sob{h}{H^k}^2}&=\Sob{h}{2}^2+\sum_{\ell=1}^k\Sob{D^\ell h}{L^2}^2.\notag
\end{align}

\vincent{From now on, we will generally drop the dependence on $\Om$ and simply write $H^k=H^k(\Om)$ and $L^p=L^p(\Om)$.} 
\vincent{We will make copious use of the Gagliardo-Nirenberg interpolation inequality, which is stated as follows: Suppose $1\leq p,q,r\leq\infty$ and $0\leq j<m$. Then for all $s>0$}
    \begin{align}\notag
 \vincent{\|D^jh\|_p\leq c\|D^mh\|_r^\alpha\|u\|_q^{1-\alpha}+c\|h\|_s,\quad \frac 1p = j+\left(\frac 1r -m\right)\alpha+\frac{1-\alpha}{q},}
\end{align}
\vincent{where $\alpha\in[j/m,1)$, for some universal constant $c$ depending on $p,q,r,s,j,m$. A particularly useful case for our analysis will be when $\al=1/2$:}
\begin{align}\label{GNS}
   \vincent{ \|D^jh\|^2_p\leq c(\|D^mh\|_r\|h\|_q+\|h\|^2_p).}
\end{align}


\vincent{We conclude this section by stating the local existence theory for \eqref{eq:main11}, \eqref{eq:main22}, which essentially follows the work of \cite{tao2013large} and \cite{zeng2019logarithmic}. Its proof will thus, be omitted here and we refer the reader to the thesis of the first author of this paper for additional details (see \cite{padithesis}).}


\begin{nthm}[Local existence and uniqueness of strong solutions]\label{thm:local:exist}
\vincent{Let $u_0,v_0\in H^2$ such that $\frac{du_0}{dx}, v_0\in H^1_0$ and $u_0\geq0$ and $\veps\geq0$. There exists $T_0>0$ and a unique solution $(u,v)$ of
 \eqref{eq:main11}, \eqref{eq:main22} satisfying $u\geq 0$ in $(x,t)\in[0,1]\times[0,T_0)$, $u(x,0)=u_0$, $v(x,0)=v_0$ such that}
    \begin{align}\label{def:strong}
    u,v\in C([0,T_0);H^2),\quad \bdy_xu,v\in L^\infty(0,T_0;H^1_0),\quad u\in L^2(0,T_0;H^3).
    \end{align}
\vincent{When $\veps>0$, then additionally $v\in L^2(0,T_0; H^3)$.}
\end{nthm}
We refer to any solution of \eqref{eq:main1}, \eqref{eq:main2} satisfying \eqref{def:strong} as a \textit{strong solution} of \eqref{eq:main1}, \eqref{eq:main2} over $[0,T)$. We say that the solution is \textit{global} if it satisfies \eqref{def:strong} for all $T>0$.

\section{Statements of Results}\label{sect:main:results}
In this section we will state our main results. \vincent{Our first main result identifies a Lyapunov structure for \eqref{eq:main1}, \eqref{eq:main2}. This will ultimately be leveraged to initiate the bootstrap from $L^2$ to $H^1$ to $H^2$. In particular, we will establish the following estimates.}

\begin{nthm}\label{thm:Lyapunov}
\vincent{Let $(u,v)$ denote the unique strong solution of \eqref{eq:main1}, \eqref{eq:main2} over $[0,T)$ corresponding to initial data $(u_0,v_0)$ satisfying}
    \begin{align}\label{def:M}
       M_0=M_0(u_0,v_0):=\intol ( u_0(x) -\ln u_0(x))dx <\infty.
    \end{align}
\vincent{Then}
    \begin{align}
   & \vincent{\frac{d}{dt}\Bigg{(} \intol (u(x,\cdotp)\ln u(x,\cdotp)-u(x,\cdotp))dx + \frac 12 \|v(\cdotp)\|^2 \Bigg{)}+\intol u_x(x,t)(\ln u(x,t))_xdx+\veps\|v_x(t)\|^2 \leq 0}\label{eq:Lyapunov1}\\
    &\vincent{\frac{d}{dt}\left(\intol (u(x,\cdotp)-\ln u(x,\cdotp))dx\right) +\frac 12 \left \|\bdy_x\ln u(t)\right\|^2+r \|u(t)-1\|^2  \leq \frac 12 \|v(t)\|^2,}\label{eq:Lyapunov2}
\end{align}
\vincent{holds for all $0\leq t<T$.} 
\end{nthm}

\vincent{Our next main result is \cref{thm:Glob:exist}, which establishes the global existence and uniqueness, as well as the time-asymptotic behavior of $(u,v)$.}


\begin{nthm}[Global existence, uniqueness, and asymptotic behavior]\label{thm:Glob:exist}
Let $u_0,v_0\in H^2$ such that $\frac{du_0}{dx}, v_0\in H^1_0$ and $u_0\geq0$ and $M_0<\infty$. Then for any $\veps\in(0,1)$ and $r>0$, there exists a unique, global strong solution $(u,v)$ of \eqref{eq:main1}, \eqref{eq:main2} corresponding to the initial data $(u_0,v_0)$ such that
    \begin{align}\label{est:H2:energy:ineq}
         \|u(t)-1\|_{H^2}^2+\|v(t)\|_{H^2}^2 +\intot \left(\|u(s)-1\|^2_{H^3}ds+\veps\|v(s)\|^2_{H^3}\right)ds\leq C_0,
    \end{align}
for all $t\geq0$, for some positive constant $C_0$ depending on $r,\veps$, but independent of $t$. Moreover, there exists $T^*>0$ such that
 \begin{align}
 \|u(t)-1\|_{H^2}^2+\|v(t)\|_{H^2}^2 \leq C_1 e^{-C_2 (t-T^*)},\label{est:decay:rate}
 \end{align}
for all $t\geq T^*$, where $C_1, C_2$ are positive constants depending on $r,\veps$, but independent of $t$.
\end{nthm}

\vincent{Our third main result is \cref{thm:diff} which establishes that the vanishing chemical diffusion limit holds.} In particular, it establishes that (strong) solutions of the system with diffusion in the chemical concentration do in fact converge to the solutions of the system without diffusion in $v$.

\begin{nthm}\label{thm:diff}(Zero chemical diffusivity limit)
Let $u_0,v_0\in H^2$ such that $\frac{du_0}{dx}, v_0\in H^1_0$ and $u_0\geq0$ and $M_0<\infty$. Given $T>0$ and $\veps\geq0$, let $(u^\veps,v^\veps)$ denote the unique strong solution of
 \eqref{eq:main1}, \eqref{eq:main2} corresponding to initial values $(u_0,v_0)$ over $[0,T]$. Then
 \begin{align}
   \sup_{t\in[0,T]}\left(\|u^\veps(t)-u^0(t)\|^2_{H^1}+\|v^\veps(t)-v^0(t)\|^2_{H^1}\right)\leq \mathcal{O}_T(\veps),\notag
 \end{align}
 where $\mathcal{O}_T(\veps)$ denotes a quantity that depends on $T$ and $\veps$ such that $\lim_{\veps\rightarrow0^+}\mathcal{O}_T(\veps)=0$.
\end{nthm}

\vincent{We will first establish the Lyapunov structure in \cref{thm:Lyapunov} in \cref{sect:Lyapunov}. Next, we will prove \cref{thm:Glob:exist} in \cref{sect:Glob:exist}. As this will rely on the local existence theory, which is done in \cite{padithesis}, we point out that it will suffice to establish the relevant {\it a priori} estimates. We will, in fact, obtain estimates that are independent of time. As usual, a standard continuation argument will then allow one to extend the solution beyond the local existence time. We will finally prove \cref{thm:diff} in \cref{sect:diff}. With \cref{thm:Glob:exist} in hand, it will suffice to establish the relevant {\it a priori} estimates, but independent of $\veps$.}

\begin{nrk}\label{rmk:constants}
\vincent{In the analysis below, we will adopt that convention that $c,C$ denote constants that may change line-to-line. Typically, $c$ will be used to denote a universal constant that does not depend on any of the model parameters or initial conditions, whereas $C$ may depend on the model parameters and initial conditions. Whenever relevant, we will express the dependencies of $C$ on the parameters as arguments of function, e.g., $C=C(r,\veps,T)$ is a constant that depends on $\veps, r$ or ``final time" $T$.} \vincent{Also, all space integrals will be denoted by $\intol$ and all time integrals will be denoted as $\int_0^t$. We will often omit the notation $dx$, $dt$.}
\end{nrk}





\section{Lyapunov structure}\label{sect:Lyapunov}

Suppose that $(u,v)$ is the unique, strong solution of \eqref{eq:main1}, \eqref{eq:main2} guaranteed by \cref{thm:local:exist} over some time interval $[0,T_0)$. To prove \cref{thm:Lyapunov}, it will be convenient to introduce the following functions: 
    \begin{align}\label{def:entropies}
        \eta(z)=z \ln(z)-z\quad\text{and}\quad f(z)= z-\ln(z).
    \end{align}
\vincent{Thus, it is equivalent to show}
\begin{align}
   & \vincent{\frac{d}{dt}\Bigg{(} \intol\eta(u(x,\cdotp)) + \frac 12 \|v(\cdotp)\|^2 \Bigg{)}+\intol \frac{u_x(t)^2}{u(t)}+\veps\|v_x(t)\|^2 \leq 0}\label{pre10}\\
    &\vincent{\frac{d}{dt}\left(\intol f(u(x,\cdotp))\right) +\frac 12 \left \|\frac{u_x(t)}{u(t)}\right\|^2+r \|u(t)-1\|^2  \leq \frac 12 \|v(t)\|^2},\label{pre20}
\end{align}
\vincent{holds for all $t\in[0,T_0)$.}

\begin{nrk}\label{V2}
\vincent{Observe that} \eqref{pre10} is equivalent to 
\[
\frac{d}{dt}\left( \intol\eta(u)-\eta(1)-\eta'(1)(u-1) + \frac 12 \|v\|^2\right)+\intol \frac{u_x^2}{u}+\veps\|v_x\|^2 \leq 0,
\]
holds for all $t\in[0,T_0)$. Notice that $\eta(z)$ is a convex function, which implies $\eta(u)-\eta(1)-\eta'(1)(u-1)\geq 0$. In particular, this gives 
    \begin{align}\label{eq:V2}
    \|v(t)\|^2\leq \|v_0\|^2=C\quad\text{and}\quad 
    \veps\intot \|v_x\|^2\leq C,
    \end{align}
for $C$ independent of $\veps$ and $t$.
\end{nrk}


\vincent{Therefore, to prove \cref{thm:Lyapunov}, we will show that \eqref{pre10} and \eqref{pre20} hold for all $t\in[0,T_0)$.}


\begin{proof}[Proof of \cref{thm:Lyapunov}]
First, we prove \eqref{pre10}. By taking $L^2$ inner product of $\ln(u)$ with \eqref{eq:main11} and of $v$ with \eqref{eq:main22} and adding it together and using integration by parts, we get
\begin{align}
\frac{d}{dt}\left( \intol\eta(u) + \frac 12 v^2\right)+\intol \frac{u_x^2}{u}+\veps\intol v_x^2 =&\intol\R  \ln(u)(1-u)u -\veps\intol v^2v_x \notag\\
&+ \left(\R \ln(u)uv+\veps v^3+\veps v v_x\right)\Big{]}_0^1\notag
\end{align}
\vincent{where $\eta(u)=u\ln(u)-u$. Notice that \eqref{eq:main22} makes the last two terms on the right-hand side equal to $0$, so that}
\begin{align}
\frac{d}{dt}\left( \intol\eta(u) + \frac 12 v^2\right)+\intol \frac{u_x^2}{u}+\veps\intol v_x^2 =\intol \R \ln(u)(1-u)u. \notag
\end{align}
\vincent{Also note} that term on the right-hand side is a non-positive term. This is clear if $u\gg1$. For $u\leq 1$, we can use the Taylor expansion to obtain, 
\begin{align}
r\ln(u)(1-u)u\leq r\frac {1}{\bar u}(u-1)(1-u)u= -r \frac {1}{\bar u}(1-u)^2u \leq0, \quad \bar u\in(0,1).\notag
\end{align}

\vincent{Similarly, to obtain \eqref{pre20}, we take the $L^2$ inner product of $\frac{u-1}{u}$ with \eqref{eq:main11} and use integration by parts to obtain}
\begin{align}
\frac{d}{dt}\left( \intol f(u)\right)+\intol \frac{u_x^2}{u^2}+\R \intol (1-u)^2 =\intol \frac{u_xv}{u}+ \left[ (u-1)\left(\frac {u_x}{u}-v\right)\right]_0^1.\notag
\end{align}
\vincent{Using \eqref{eq:main11} and \eqref{eq:main22} we obtain}
\begin{align}
\frac{d}{dt}\left( \intol f(u)\right)+\intol \frac{u_x^2}{u^2}+\R \intol (1-u)^2 =\intol \frac{u_xv}{u}.\notag
\end{align}
Notice that the function $f(u)$ is non-negative and convex. \vincent{Now} we can do the following estimate using the Cauchy-\vincent{Schwarz} inequality
\begin{align}
\frac{d}{dt}\left( \intol f(u)\right)+\intol \frac{u_x^2}{u^2}+\R \intol (1-u)^2 =\intol \frac{u_xv}{u} \leq \frac 12 \|\frac{u_x}{u}\|^2 +\frac 12 \|v\|^2.\notag
\end{align}
Thus,
\begin{align}\label{eqn3}
\frac{d}{dt}\left( \intol f(u)\right)+\frac 12 \|\frac{u_x}{u}\|^2+\R \|1-u\|^2 \leq \frac 12 \|v\|^2,
\end{align} 
\vincent{as desired}.
\end{proof}

\section{Global existence and uniqueness: Proof of \cref{thm:Glob:exist}}\label{sect:Glob:exist}

\vincent{We develop a standard bootstrap procedure. In particular, we obtain $L^2$ estimates, then $H^1$ estimates, then $H^2$ estimates. In the analysis we perform below, it will be convenient to introduce the shifted variable, $\tu=u-1$. Then \eqref{eq:main11} and \eqref{eq:main22} become}
\begin{align}\label{shift:eqn:u}
\tu_t +(\tu v)_x+v_x&=\tu_{xx}-\R (\tu +1)\tu\\
v_t+\tu_x&=\veps v_{xx}+\veps(v^2)_x\label{shift:eqn:v}
\end{align}
with initial and boundary conditions given by
\begin{align}\label{ic1}
&\tu_0> -1,\quad \tu_0,v_0\in H^2(\Omega);\\
&\partial_x\tu|_{\partial\Omega}=0,\quad  v|_{\partial\Omega}=0\label{NBC1}.
\end{align}

With this notation, we may rewrite \cref{thm:Lyapunov}:
\begin{align}
    \frac{d}{dt}\Bigg{(} \intol\eta(\tu +1) + \frac 12 \|v\|^2 \Bigg{)}+\intol \frac{\tu_x^2}{\tu +1}+\veps\|v_x\|^2 &\leq 0\label{pre1}\\
    \frac{d}{dt}\left(\intol f(\tu +1)\right) +\frac 12 \left \|\frac{\tu_x}{\tu+1}\right\|^2+r \|\tu\|^2  &\leq \frac 12 \|v\|^2.\label{pre2}
\end{align}

\begin{nrk}\label{rk:extrabc}
Notice that by \eqref{NBC1}, we obtain the following compatibility condition
    \begin{align}\notag
       \vincent{ \tu_{xt}|_{\partial \Omega}=0=v_t|_{\partial \Omega}.}
    \end{align}
Also, by \eqref{shift:eqn:v} and taking $\partial_x$ of \eqref{shift:eqn:u}, we obtain the additional compatibility condition
    \begin{align}\notag
   \vincent{ v_{xx}|_{\partial \Omega}=0=\tu_{xxx}|_{\partial\Omega}.}
    \end{align}
These compatibility conditions play a crucial role in the analysis below. In contrast, such additional conditions are not required to treat the case of the unbounded domain, $\Om=\mathbb{R}$, as in \cite{zeng2019logarithmic}. 
\end{nrk}

\subsection{\vincent{$L^2$ estimates}}
\begin{nlem}[$L^2$ Estimate]\label{lem:L2} Let $\tu, v$ be strong solutions of equations \eqref{shift:eqn:u}, \eqref{shift:eqn:v} with \eqref{ic1}, \eqref{NBC1}. Then
\begin{align}\label{est:L2}
    \|\tu(t)\|^2+\|v(t)\|^2 +\int_0^t(\|\tu_{x}\|^2+\veps\|v_{x}\|^2)+\frac{c}{\R}\intot \left \|\frac{\tu_x}{\tu+1}\right\|^2+\frac12\intot \|\tu\|^2\leq \frac{C}{\veps r},
\end{align}
\vincent{for all $t\geq0$, where $c,C>0$ are independent of $\veps,r,t$}.
\end{nlem}

\begin{proof}[Proof of Lemma \ref{lem:L2}]
Taking the $L^2$ inner products of \eqref{shift:eqn:u} with  $u$ and \eqref{shift:eqn:v} with  $v$ and adding the results, we get
\begin{align}\label{L2preest}
 \frac 12\frac{d}{dt} \left(  \|\tu\|^2 + \|v\|^2 \right)+ \|\tu_x\|^2 +\veps\|v\|^2_x + r\intol \tu^2 (\tu+1)= \intol\tu_x\tu v.
\end{align}
\vincent{We estimate the right-hand side with the Cauchy-Schwarz inequality, H\"older's inequality, and Young's inequality to get}
\begin{align}
 \intol\tu_x\tu v\leq \|\tu_x\| \|\tu v\| \leq & \frac 14 \|\tu_x\|^2+ \|\tu v\|^2\notag\\
 \leq &\frac 14 \|\tu_x\|^2+ \|\tu\|_{\infty}^2\|v\|^2\notag\\
 \leq &\frac 14 \|\tu_x\|^2+c \|\tu\|_{\infty}^2.\notag
\end{align}
Notice that in the last step we have used \vincent{\eqref{eq:V2} from \cref{V2}}. Now, by \eqref{GNS} we have
\begin{align}
c \|\tu\|_\infty^2\leq& c(\|\tu_x\|\|\tu\|+\|\tu\|^2)\notag
\leq \frac14\|\tu_x\|^2+c\|\tu\|^2, \notag
\end{align}
which implies
\begin{align}\notag
 \intol\tu_x\tu v\leq \frac12\|\tu_x\|^2 +c\|\tu\|^2.
\end{align}
\vincent{Returning to \eqref{L2preest}, we obtain}
\begin{align}\label{L2}
\frac{d}{dt} \left( \|\tu\|^2 + \|v\|^2 \right)+ \|\tu_x\|^2 +2\veps\| v_x\|^2 + 2r\intol \tu^2 (\tu+1)\leq c\|\tu\|^2.
\end{align}

\vincent{Now we multiply $\frac{2c}{\R}$ with \eqref{pre2} and add the result to \eqref{L2} to obtain}
\begin{align}
    \frac{d}{dt}\left( \frac{2c}{\R}\intol f(\tu+1) + \|\tu\|^2+\|v\|^2\notag \right)+\|\tu_x\|^2 +2\veps\|v_x\|^2 +&\frac{c}{\R}\left \|\frac{\tu_x}{\tu+1}\right\|^2+c \|\tu\|^2 \leq  \frac{c}{\R}\|v\|^2.\notag
\end{align}
\vincent{The} Poincar\'e inequality gives
\begin{align}
    \frac{d}{dt}\left( \frac{2c}{\R}\intol f(\tu+1) + \|\tu\|^2+\|v\|^2 \right)+\|\tu_x\|^2 +2\veps\|v_x\|^2 +&\frac{c}{\R}\left \|\frac{\tu_x}{\tu+1}\right\|^2+c \|\tu\|^2 \leq  \frac{c}{\R}\|v_x\|^2.\notag
\end{align}
 Integrating in time, using the positivity of $f$, \eqref{pre1}, and \vincent{applying \cref{thm:Lyapunov}},  we obtain 
\begin{align}
   \|\tu(t)\|^2+\|v(t)\|^2 +\intot\|\tu_x\|^2 +2\veps\intot\|v_x\|^2 +\frac{c}{\R}\intot \left \|\frac{\tu_x}{\tu+1}\right\|^2+c\intot \|\tu\|^2\leq \frac{C}{\veps r}
   ,\notag
\end{align}
\vincent{which implies \eqref{est:L2}, as claimed.}
\end{proof}
\subsection{\vincent{$H^1$ estimates}}
\begin{nlem}[$H^1$ Estimate]\label{lem:H1} Let $\tu, v$ be strong solutions of the equations \eqref{shift:eqn:u}, \eqref{shift:eqn:v} with \eqref{ic1} and \eqref{NBC1}. Then
\begin{align}\label{est:H1}
    \|\tu_x(t)\|^2+\|v_x(t)\|^2 +\int_0^t(\|\tu_{xx}\|^2+\veps\|v_{xx}\|^2)\leq C(r,\veps),
\end{align}
where $C(r,\veps)$ depends on $\veps, r$, but is independent of $t$. In particular, $C(r,\veps)=C\exp\left({c\frac{r^{-1}\vee r}{\veps}}\right)$, for some $c,C>0$ independent of $\veps, r,t$.
\end{nlem}

\begin{proof}[Proof of Lemma \ref{lem:H1}]
First \vincent{we will} estimate $\|v_x\|^2$. \vincent{To do so, we apply $\partial_x$ of \eqref{shift:eqn:v} and add the result to \eqref{shift:eqn:u}, which gives}
\begin{align}
v_{tx}=\veps v_{xxx}+\veps (v^2)_{xx}-\tu_t-(\tu v)_x-v_x-\R (\tu+1)\tu. \notag
\end{align}
\vincent{Taking the} $L^2$ inner product with $v_x$ and \vincent{integrating} by parts \vincent{gives}
\begin{align}\label{eq:H1:balance:v}
\frac{d}{dt}\left( \frac12\|v_x\|^2\right)+\veps\|v_{xx}\|^2+\|v_x\|^2=-\veps \intol v_{xx}(v^2)_{x}-\intol v_x\tu_t-\intol v_x(\tu v)_x-\R\intol v_x(\tu+1)\tu.
\end{align}
Using \vincent{the identity} $\frac{d}{dt}\left(\intol \tu v_x\right)=\intol \tu_t v_x +\intol \tu v_{xt}$, \vincent{we can rewrite \eqref{eq:H1:balance:v} as}
\begin{align}\label{eq:H1:balance:v2}
    \begin{split}
\frac{d}{dt}&\left( \frac12\|v_x\|^2+\intol \tu v_x\right)+\veps\|v_{xx}\|^2+\|v_x\|^2\\
&=
-\veps \intol v_{xx}(v^2)_{x}+\intol v_{xt}\tu-\intol v_x(\tu v)_x-\R\intol v_x(\tu+1)\tu.
    \end{split}
\end{align}
\vincent{From the equation, observe that $v_{xt}=-\tu_{xx}+\veps v_{xxx}+\veps(v^2)_{xx}$. Using this identity in \eqref{eq:H1:balance:v2}, we get}
\begin{align}
\frac{d}{dt}\left( \frac12\|v_x\|^2+\intol \tu v_x\right)&+\veps\|v_{xx}\|^2+\|v_x\|^2 \notag\\
=-\veps \intol v_{xx}(v^2)_{x}-&\intol\tu \tu_{xx}+\veps\intol \tu v_{xxx}+\veps\intol(v^2)_{xx}\tu-\intol v_x(\tu v)_x-\R\intol v_x(\tu+1)\tu. \notag
\end{align}
\vincent{Integrating by parts again and re-grouping terms}, we obtain
\begin{align}
\frac{d}{dt}&\left( \frac12\|v_x\|^2+\intol \tu v_x\right)+\veps\|v_{xx}\|^2+\|v_x\|^2 \notag\\
&=\vincent{\| \tu_{x}\|^2-2\veps \intol v_{xx}v v_{x}-\veps\intol \tu_x (v_{xx}+2\tu v_xv)-\intol v_x(\tu_x v+\tu v_x)-\R\intol v_x(\tu+1)\tu} \notag\\
&=\vincent{\|\tu_x\|^2+V_1+V_2+V_3+V_4.}\label{eq:H1:regroup}
\end{align}
We can estimate \vincent{each of the} terms on the right-hand side using \vincent{H\"older's inequality}, Young's inequality, \cref{lem:L2}, and \eqref{GNS}. \vincent{For $V_1$, we estimate}
\begin{align}
|V_1|
&\leq \frac \veps 8\|v_{xx}\|^2+c\veps\|v\|^2\|v_x\|_\infty^2 \notag\\
&\leq \frac \veps 8\|v_{xx}\|^2+c\veps (\|v_{xx}\|\|v_x\|+\|v_x\|^2) \notag\\
&\leq \frac \veps 4\|v_{xx}\|^2+c\veps \|v_x\|^2.\notag
\end{align}
\vincent{For $V_2$, we estimate}
\begin{align}
|V_2|
\leq &  \frac\veps 8\|v_{xx}\|^2 +c\veps\|\tu_x\|^2 +\veps \|v\|^2\|v_x\|^2_\infty+\veps\|\tu_x\|^2\notag \\
\leq &  \frac\veps 8\|v_{xx}\|^2 +c\veps\|\tu_x\|^2 +c\veps(\|v_x\|\|v_{xx}\|+\|v_x\|^2)\notag \\
\leq &  \frac\veps 4\|v_{xx}\|^2 +c\veps\|\tu_x\|^2 +c\veps\|v_x\|^2.\notag 
\end{align}
\vincent{For $V_3$, we estimate}
\begin{align}
|V_3|
\leq &\frac 14 \|v_x\|^2 + (\|\tu\|^2_\infty \|v_x\|^2+\|v\|_\infty^2\|\tu_x\|^2)\notag\\
\leq &\frac 14 \|v_x\|^2 +c (\|\tu_x\| \|\tu\|\|v_x\|^2+\|\tu\|^2\|v_x\|^2+\|v_x\|\|v\|\|\tu_x\|^2+\|v\|^2\|\tu_x\|^2)\notag\\
\leq &\frac 14 \|v_x\|^2 + c (\|\tu_x\|^2+\|\tu\|^2)\|v_x\|^2+c\|\tu_x\|^2.\notag
\end{align}
\vincent{For $V_4$, we estimate}
\begin{align}
|V_4|
\leq&  c\R^2\|\tu\|^2+\frac 14 \|v_x\|^2+\R \|\tu\|_\infty^2\|v_x\| \notag\\
\leq &c\R^2\|\tu\|^2+\frac 14 \|v_x\|^2+ \R c(\|\tu_x\|\|\tu\|+\|\tu\|^2)\|v_x\|\notag\\
\leq &c\R^2\|\tu\|^2+\frac 14 \|v_x\|^2+\R^2c\|\tu\|^2+c\|\tu_x\|^2\|v_x\|^2+\R c\|\tu\|^2\|v_x\| \notag\\
\leq &c\R^2\|\tu\|^2+\frac 14 \|v_x\|^2+c(\|\tu_x\|^2+r\|\tu\|^2)\|v_x\|^2.\notag
\end{align}

\vincent{Combining $V_1$--$V_4$ in \eqref{eq:H1:regroup} we get}
\begin{align}
\frac{d}{dt}\left( \frac12\|v_x\|^2+\intol \tu v_x\right)&+\frac\veps2\|v_{xx}\|^2+\frac 12\|v_x\|^2 \notag\\
\leq c&(\|\tu_x\|^2+ (1+r)\|\tu\|^2)\|v_x\|^2+c(\veps\|v_x\|^2+\|\tu_x\|^2+\R^2 \|\tu\|^2).\notag
\end{align}
Adding \eqref{L2} we obtain
\begin{align}
\frac{d}{dt}&\left( \|\frac 12 v_x+\tu\|^2+\|v\|^2+\frac 14\| v_x\|^2 \right) +\frac\veps2\|v_{xx}\|^2+\frac 12\|v_x\|^2+\|\tu_x\|^2+\veps\|v_x\|^2 \notag\\
\leq & c(\|\tu_x\|^2+ (1+r)\|\tu\|^2)\frac 14\|v_x\|^2 + c(\veps\|v_x\|^2+\|\tu_x\|^2+(1+\R^2) \|\tu\|^2).\notag
\end{align}
\vincent{By \cref{lem:L2}, we have
$\int_0^t\|\tu\|^2, \int_0^t\|\tu_x\|^2\leq C(\veps r)^{-1}$, $\intot\|v_x\|^2\leq C(\veps^2\R)^{-1}$. Thus, by Gr\"onwall's inequality to obtain}
\begin{align}\label{vH1}
\frac 14 \|v_x\|^2+\frac\veps2\int_0^t\|v_{xx}\|^2+\frac 12\int_0^t\|v_x\|^2+\int_0^t\|\tu_x\|^2+\veps\int_0^t\|v_x\|^2\leq \frac{C}{\veps}(r^{-1}\vee r)^2e^{\frac {c}{\veps r}} \le Ce^{c\frac{r^{-1}\vee r}{\veps}},
\end{align}
where $a\vee b=\max\{a,b\}$. Indeed, observe that $a^{-1}\vee a\geq1$, for all $a>0$.


\vincent{We will now estimate $\|\tu_x\|^2$. To do so, we take $\partial_x$ of \eqref{shift:eqn:u} and \eqref{shift:eqn:v}, then take the $L^2$ inner product with $\tu_x$ and $v_x$, respectively, and sum the results. After integrating by parts we get}
\begin{align}\label{eq:H1:balance:u}
\frac 12 \frac{d}{dt}(\|\tu_x\|^2+\|v_x\|^2)+\|\tu_{xx}\|^2+\veps\|v_{xx}\|^2&=\intol (\tu_x v+\tu v_x)\tu_{xx}+\R\intol \tu(\tu+1)\tu_{xx}-2\veps\intol vv_xv_{xx}\notag\\
&=U_1+U_2+U_3.
\end{align}
\vincent{Using H\"older's inequality, \eqref{GNS}, and Young's inequality, we get}
\begin{align}
U_1
\leq& \frac{1}{8} \|\tu_{xx}\|^2+c\|v\|_\infty^2\|\tu_x\|^2+C\|\tu\|_\infty^2\|v_x\|^2\notag\\
\leq& \frac18 \|\tu_{xx}\|^2 +c(\|v_x\|\|v\|+\|v\|^2)\|\tu_x\|^2+c(\|\tu_x\|\|\tu\|+\|\tu\|^2)\|v_x\|^2\notag\\
\leq& \frac18 \|\tu_{xx}\|^2 +c(\|v_x\|^2+\|v\|^2)\|\tu_x\|^2+c(\|\tu_x\|^2+\|\tu\|^2)\|v_x\|^2.\notag
\end{align}
\vincent{Using the Cauchy-Schwarz inequality, \eqref{GNS}, and \cref{lem:L2}, we get}
\begin{align}
U_2
&\leq r \|\tu\|_\infty^2\|\tu_{xx}\|^2+r\intol \vincent{|\tu||\tu_{xx}|}\notag\\
&\leq cr(\|\tu_x\|\|\tu\|+\|\tu\|^2)\|\tu_{xx}\|+cr^2\|\tu\|^2+\frac18\|\tu_{xx}\|^2\notag\\
&\leq cr^2\|\tu\|^2\|\tu_x\|^2+\frac38\|\tu_{xx}\|^2 +cr^2\|\tu\|^4+cr^2\|\tu\|^2\notag\\
&\leq cr^2\|\tu\|^2\|\tu_x\|^2+\frac38\|\tu_{xx}\|^2 +\vincent{C\left(\frac{r^{-1}\vee r}{\veps}\right)^2\|\tu\|^2}\notag.
\end{align}
\vincent{By H\"older's inequality, \eqref{GNS}, and Young's inequality, we get}
\begin{align}
U_3
\leq & \veps(\frac 12 \|v_{xx}\|^2+c\|v\|^2_\infty\|v_x\|^2)\notag\\
\leq & \frac \veps 2 \|v_{xx}\|^2+c\veps (\|v_x\|\|v\|+\|v\|^2)\|v_x\|^2\notag\\
\leq & \frac \veps 2 \|v_{xx}\|^2+c\veps \|v_x\|^2\|v_x\|^2 +c\veps \|v_x\|^2.\notag
\end{align}
{\vincent{Combining $U_1$--$U_3$ in \eqref{eq:H1:balance:u}, we get}}
\begin{align}\label{H1}
\frac 12 &\frac{d}{dt}(\|\tu_x\|^2+\|v_x\|^2)+\frac 12\|\tu_{xx}\|^2+\frac\veps 2\|v_{xx}\|^2\\
\leq & c(\|v_x\|^2+\|v\|^2+r^2\|\tu\|^2+\|\tu_x\|^2+\veps\|v_x\|^2)(\|\tu_x\|^2+\|v_x\|^2)+C\vincent{\left(\frac{r^{-1}\vee r}{\veps}\right)^2}\|\tu\|^2+C\veps\|v_x\|^2.\notag
\end{align}
\vincent{Observe that from \cref{lem:L2}, the Poincar\'e inequality, and \eqref{vH1}, we have $\int_0^t \|v\|^2\leq\int_0^t \|v_x\|^2\leq C(\veps^2r)^{-1}$ and
$\int_0^t\|\tu\|^2,\int_0^t \|\tu_x\|^2\leq C(\veps r)^{-1}$. Thus, by Gr\"onwall's inequality, we deduce}
\begin{align}
\|\tu_x\|^2+\|v_x\|^2 +\frac 12\int_0^t\|\tu_{xx}\|^2+\frac\veps 2\int_0^t\|v_{xx}\|^2 \leq \vincent{Ce^{c\frac{r^{-1}\vee r}{\veps}}},\notag
\end{align}
\vincent{which implies \eqref{est:H1}, as desired.}
\end{proof}

\subsection{\vincent{$H^2$ estimates}}
\begin{nlem}[$H^2$ Estimate]\label{lem:H2} Let $\tu, v$ be solutions of the equations \eqref{shift:eqn:u}, \eqref{shift:eqn:v} with \eqref{ic1} and \eqref{NBC1}. Then, it follows that 
\begin{align}\label{est:H2}
    \|\tu_{xx}(t)\|^2+\|v_{xx}(t)\|^2+\intot (\|\tu_{xxx}\|^2+\veps \|v_{xxx}\|^2)\leq C(r,\veps),
\end{align}
for all $t\geq0$, where $C=C(r,\veps)$ is independent of $t$. In particular, $C(r,\veps)=C\exp\left({c\frac{r^{-1}\vee r}{\veps^2}}\right)$, for some constants $c,C>0$, independent of $\veps,r$.
\end{nlem}

\begin{proof}[Proof of \cref{lem:H2}]

\vincent{We apply $\partial_{xx}$ to \eqref{shift:eqn:u} and \eqref{shift:eqn:v}, then take the $L^2$-inner product with $\tu_{xx}$ and $v_{xx}$, respectively, and sum the results. After integrating by parts and using \cref{rk:extrabc}, we obtain}
\begin{align}
\frac 12\frac{d}{dt}\left( \|\tu_{xx}\|^2+\|v_{xx}\|^2\right)+&\|\tu_{xxx}\|^2+\veps\|v_{xxx}\|^2\notag\\
=&\intol (\tu v)_{xx}\tu_{xxx}+\R \intol((\tu+1)\tu)_{x}\tu_{xxx}-\veps\intol(v^2)_{xx}v_{xxx}\notag\\
=&\intol (\tu v_{xx}+2\tu_x v_x+v \tu_{xx})\tu_{xxx}+\R\intol (2\tu \tu_x +\tu_x)\tu_{xxx}-2\veps\intol(vv_{xx}+v^2_x)v_{xxx}\notag\\
=& W_1+W_2+W_3\label{eq:H2:balance:u}.
\end{align}
\vincent{We can treat each term by making use of integration by parts, \cref{rk:extrabc}, as well as H\"older's inequality, \eqref{GNS}, and \cref{lem:H1}. For $W_1$, we estimate}
\begin{align}
W_1
&\leq \frac 14\|\tu_{xxx}\|^2+ c\|\tu\|^2_\infty\|v_{xx}\|^2+c\|v\|_\infty^2\|\tu_{xx}\|^2+c\|v_x\|_\infty^2\|\tu_x\|^2\notag\\
&\leq \frac 14\|\tu_{xxx}\|^2+ c(\|\tu\|^2+\|\tu_x\|^2)\|v_{xx}\|^2+c(\|v\|^2+\|v_x\|^2)\|\tu_{xx}\|^2+c(\|v_x\|^2+\|v_{xx}\|^2)\|\tu_x\|^2\notag\\
&\leq \frac 14\|\tu_{xxx}\|^2+ c(\|\tu\|^2+\|\tu_x\|^2)\|v_{xx}\|^2+c(\|v\|^2+\|v_x\|^2)\|\tu_{xx}\|^2+c\|v_x\|^4+c\|\tu_x\|^4.\notag
\end{align}
\vincent{For $W_2$, we estimate}
\begin{align}
    W_2
    &\leq \frac14 \|\tu_{xxx}\|^2+c\R^2\|\tu\|_\infty^2\|\tu_x\|^2+c\R^2\|\tu_x\|^2\notag\\
    &\leq \frac14 \|\tu_{xxx}\|^2+cr^2(1+\|\tu\|^4+\|\tu_x\|^4).\notag
\end{align}
\vincent{For $W_3$, we use the fact that $\veps\leq1$ and estimate}
\begin{align}
    W_3
    &\leq \frac\veps 2\|v_{xxx}\|^2+c\veps\|v\|_\infty^2\|v_{xx}\|^2+c\veps\|v_x\|_\infty^2\|v_x\|^2\notag\\
    &\leq \frac\veps 2\|v_{xxx}\|^2+c\veps(\|v\|^2+\|v_x\|^2)\|v_{xx}\|^2+c\veps(\|v_x\|^2+\|v_{xx}\|^2)\|v_x\|^2\notag\\
     &\leq \frac\veps 2\|v_{xxx}\|^2+c\veps(\|v\|^2+\|v_x\|^2)\|v_{xx}\|^2+c\veps\|v_x\|^4.\notag
\end{align}
\vincent{Combining $W_1$--$W_3$ in \eqref{eq:H2:balance:u}}, we obtain
\begin{align}\label{est:H2:ineq}
   \frac12\frac{d}{dt} (\|\tu_{xx}\|^2+\|v_{xx}\|^2)+&\frac12\|\tu_{xxx}\|^2+\frac\veps2\|v_{xxx}\|^2
    \leq M(t)(\|\tu_{xx}\|^2+\|v_{xx}\|^2)+N(t),
\end{align}
where 
\begin{align}
    M(t)&:=c(\|\tu\|^2+\|\tu_x\|^2+\|v\|^2+\|v_x\|^2),\notag\\
    N(t)&:=c(r^{-1}\vee r)^2(\|v_x\|^4+ \|\tu_x\|^4).\notag
\end{align}
\vincent{Note that we have applied the assumption that $\veps\leq1$. Observe from \cref{lem:L2} and \cref{lem:H1} that $\intot M(t)\leq C(\veps^2r)^{-1}$, and  $\intot N(t) \leq C(\veps r)^{-1}e^{2c\frac{r^{-1}\vee r}{\veps}}$}. We can apply Gr\"onwall's inequality to obtain
\begin{align}
    \|\tu_{xx}\|^2+\|v_{xx}\|^2 +\intot (\|\tu_{xxx}\|^2+\veps\|v_{xxx}\|^2)\leq Ce^{c\frac{r^{-1}\vee r}{\veps^2}},\notag
\end{align}
which implies \eqref{est:H2}.
\end{proof}

As an immediate corollary of the above estimates, we deduce the desired uniform-in-time bounds in $H^2$ for $(\tu, v)$.

\begin{ncor}\label{cor:H2:complete}
Let $\tu, v$ be strong solutions of the equations \eqref{shift:eqn:u}, \eqref{shift:eqn:v} with \eqref{ic1} and \eqref{NBC1}. Then
\begin{align}\label{est:H2:complete}
    \|\tu(t)\|_{H^2}^2+\|v(t)\|_{H^2}^2 +\int_0^t(\|\tu(s)\|_{H^3}^2+\veps\|v(s)\|_{H^3}^2)\leq C(r,\veps),
\end{align}
where $C(r,\veps)$, depends on $\veps, r$, but is independent of $t$. In particular, $C(r,\veps)=C\exp\left({c\frac{r^{-1}\vee r}{\veps^2}}\right)$, for some $c,C>0$ independent of $\veps, r,t$.
\end{ncor}
\subsection{Asymptotic behavior} 
\vincent{Now we will establish the asymptotic behavior asserted in \cref{thm:Glob:exist}}. We first prove that $u$ converges to its carrying capacity (normalized \vincent{to} 1) in the \vincent{uniform topology}. We will \vincent{then use this to prove} \cref{bounds} below; \vincent{this will be crucial in obtaining} the asymptotic behavior of the solutions.

\vincent{Let us first recall that} $f(z)=z-\ln z\geq 0$ and $\eta(z)=z\ln z -z $, \vincent{both of which are} convex. \vincent{In particular, recall that from  \eqref{pre2} and Poincar\'e's inequality, we have}
\begin{align}\label{pre22}
    \frac{d}{dt}\left(\intol f(\tu+1)\right) +\frac 12 \left \|\frac{\tu_x}{\tu+1}\right\|^2+\R \|\tu\|^2  \leq \frac 12 \|v\|^2\leq c\|v_x\|^2,
\end{align}
\vincent{and also from \eqref{pre1} that}
\begin{align}\label{pre11}
    \frac{d}{dt}\Bigg{(} \intol\eta(\tu+1) +& \frac 12 \|v\|^2 \Bigg{)}+\intol \frac{\tu_x^2}{\tu+1 }+\veps\|v_x\|^2 \leq 0.
\end{align}
\begin{nlem} \label{smallu}
 Under the assumptions of \cref{thm:Glob:exist}, we have
    \begin{align}\notag
    \lim_{t\rightarrow\infty}(\|\tu(t)\|_\infty+\|v(t)\|_\infty)=0.
    \end{align}
\end{nlem}
\begin{proof}[ Proof of \cref{smallu}]
Multiplying \eqref{pre22} by $\frac{2c}{\veps}$, and adding \eqref{pre1}, we get
\begin{align}
\frac{d}{dt}\left(\intol f(\tu +1)+\frac{2c}{\veps}\intol \eta(\tu +1)+\frac{c}{\veps}\|v\|^2 \right) +\frac 12 \left \|\frac{\tu_x}{\tu+1}\right\|^2+\R \|\tu\|^2 +\frac{2c}{\veps}\intol \frac{\tu^2_x}{\tu +1} +c\|v_x\|^2 \leq 0. \notag
\end{align}
\vincent{It follows that $\int_0^t (\|\tu\|^2+\|v_x\|^2)\leq C(r^{-1}\vee r)$, where $a\vee b:=\max\{a,b\}$. Using the Poincar\'e inequality, we then deduce that}
\begin{align}\notag
\intot (\|\tu\|^2+\|v\|^2)\leq C(r^{-1}\vee r),\quad \text{for all } t \in [0,\infty).
\end{align}
From \eqref{L2preest}
\begin{align}
    \left| \frac{d}{dt} (\|\tu\|^2+\|v\|^2)\right| &\leq 2\left( \|\tu_x\|^2+\veps\|v_x\|^2+\R\intol \tu^2(\tu+1)+\intol |\tu_x||\tu v|\right)\notag \\
    &\leq 2 \left( \|\tu_x\|^2 +\veps\|v_x\|^2 +\R\|\tu\|^2(\|\tu\|_\infty +1)+ \|\tu_x\|^2+\|\tu\|^2_\infty\|v\|^2\right)\notag\\
    &\leq \vincent{C(r,\veps)(\|\tu\|^2_{H^1} +\|v\|^2_{H^1})},\notag
\end{align}
\vincent{where $C(r,\veps)$ is independent of $t$.
This gives that}
\begin{align}\label{eq:F:W11}
\intot \left| \frac{d}{dt} (\|\tu\|^2+\|v\|^2)\right|\leq C(r,\veps)\intot(\|\tu\|^2_{H^1} +\|v\|^2_{H^1})\leq C(\R,\veps),
\end{align}
\vincent{holds for all $t\geq0$.}
\vincent{If we let $K(t):=\|\tu(t)\|^2+\|v(t)\|^2$, then \eqref{eq:F:W11} shows that $K(t)\in W^{1,1}(0,\infty)$. In particular, $K(t)\to 0$ as $t\to \infty$.}

\vincent{Finally, by making use of \eqref{GNS} and \cref{lem:H1}, we obtain}
\begin{align}
\|\tu(t)\|_\infty +\|v(t)\|_\infty &\leq \|\tu(t)\|^{1/2}\|\tu_x(t)\|^{1/2} +\|\tu(t)\|+\|v(t)\|^{1/2}\|v_x(t)\|^{1/2}+\|v(t)\|\notag\\
&\leq  C(r,\veps)K(t)^{1/4}+K(t)^{1/2},\notag
\end{align}
which implies $\|\tu(t)\|_\infty +\|v(t)\|_\infty \to 0$ as $t\to \infty$, \vincent{as desired}.
\end{proof}
\begin{nlem}\label{bounds} Under the assumptions of \cref{thm:Glob:exist}, there exists $T^*>0$ \vincent{such that} for all $t\geq T^* $
\begin{align}
   \frac29 \|u(t)-1\|^2 &\leq \intol f(u(t))-f(1) \leq 2\| u(t)-1\|^2,\label{f:equiv}\\
     \frac13 \|u(t)-1\|^2 &\leq \intol \eta(u(t))-\eta(1) \leq  \| u(t)-1\|^2.\label{g:equiv}
\end{align}

\end{nlem}

\begin{proof}[Proof of \cref{bounds}]
\vincent{Observe that} $f(z)-f(1)=z-\ln z-1$. \vincent{Now for each $\al>0$, define} 
    \begin{align}\label{def:g}
        g_\al(z):= z-\ln z -1-\alpha (z-1)^2.
    \end{align}
It is easy to see that 
\begin{align}\label{props:g}
\begin{split}
g_\al(1)&=0\\
g_\al^\prime(z)&=1-\frac 1z -\frac \alpha 2(z-1), \quad g'(1)=0\\
g_\al^{\prime\prime}(z)&=\frac 1{z^2}-\frac\alpha 2.
\end{split}
\end{align}
\vincent{Now observe that $g_\al(z)\geq 0$, whenever $
g_\al^{\prime\prime}(z)\geq 0$, and $g_\al(x)\leq0$, whenever $g^{\prime\prime}(z)\leq 0$. Thus, $g_\al(x)\geq0$ whenever $\alpha\leq \frac 1{2z^2}$ and $g_\al(z)\leq0$ whenever $\alpha\geq \frac 1{2z^2}.$}

Notice that as a direct implication of Lemma \ref{smallu}, there exists $T^*>0$ \vincent{such that $\|\tu(t)\|_\infty \leq \frac 12$ for all $t\geq T^* $. Since $\tu=u-1$, this} implies that for all $t\geq T^*$
\begin{align}\notag
 \frac12 \leq \|u(t)\|_\infty\leq \frac32,
\end{align}
\vincent{or equivalently that}
    \begin{align}\label{cond:u:Linfty:squared}
        \frac 29 \leq \frac 1 {2\|u(t)\|_\infty^2}\leq 2,
    \end{align}
\vincent{for all $t\geq T^*$. In particular, this implies that $g_{2}(u(t))\leq0\leq g_{2/9}(u(t))$, for all $t\geq T^*$. This proves \eqref{f:equiv}.}

\vincent{To prove \eqref{g:equiv}, first observe that} $\eta(z)-\eta(1)=z\ln z -z+1$. \vincent{Now define} 
    \begin{align}\label{def:h}
        h_\beta(z):=z\ln z -z+1-\beta(z-1)^2
    \end{align}
\vincent{It is easy to see that}
\begin{align}\label{props:h}
\begin{split}
h_\beta(\vincent{1})&=0\\
h_\beta^\prime(z)&=\ln z -2\beta(z-1),\quad \vincent{h_\beta^\prime(1)=0}\\
h_\beta^{\prime\prime}(z)&=\frac 1z- 2\beta.
\end{split}
\end{align}
\vincent{As before, we have obtain $h_\beta(z)\geq 0$, whenever
$h^{\prime\prime}(z)\geq 0$, and $h_\beta(z)\leq0$, whenever $h_\beta^{\prime\prime}(z)\leq0$. Hence, $h_\beta(z)\geq0$, whenever $\beta\leq \frac 1{2z}
$, and $h_\beta(z)\leq0$, whenever $\alpha\geq \frac 1{2z}$.}

\vincent{Thus, for $t\geq T^*$, we have}
\begin{align}
 \frac12 \leq \|u(t)\|_\infty\leq \frac32\notag,
\end{align}
\vincent{which is equivalent to}
    \begin{align}\label{cond:u:Linfty}
    \frac 13 \leq \frac 1 {2\|u(t)\|_\infty}\leq 1.
    \end{align}
\vincent{This shows that for all $t\geq T^*$, we have $h_1(u(t))\leq0\leq h_{1/3}(u(t))$, from which we deduce \eqref{g:equiv}}.
\end{proof}
Now we are ready to establish the asymptotic estimate \eqref{est:decay:rate} from \cref{thm:Glob:exist}.
\begin{proof}[Proof of \eqref{est:decay:rate} from \cref{thm:Glob:exist}]

\vincent{For convenience, let us rewrite \eqref{pre22} and \eqref{pre11} here} as 

\begin{align}
    \frac{d}{dt}\left(\intol f(\tu+1) -f(1)\right) +\frac 12 \left \|\frac{\tu_x}{\tu+1}\right\|^2+\R \|\tu\|^2  &\leq \frac 12 \|v\|^2\leq c\|v_x\|^2 \label{pre222}\\
    \frac{d}{dt}\Bigg{(} \intol (\eta(\tu+1) -\eta(1))+ \frac 12 \|v\|^2 \Bigg{)}+\intol \frac{\tu_x^2}{\tu+1 }+\veps\|v_x\|^2 &\leq 0\label{pre111}.
    \end{align}
\vincent{Multiplying \eqref{pre222} by $\frac{1}{\veps}$, then adding the result to \eqref{pre111}}, we get
\begin{align}\label{extra}
    \frac{d}{dt}&\left(\frac1\veps\intol (\eta(\tu+1)-\eta(1)) +\frac 1{2\veps}\|v\|^2+ \intol f(\tu+1)-f(1) \right)\notag\\
    &+\frac1\veps\intol \frac{(\tu_x)^2}{\tu+1}+\left\|\frac{\tu_x}{\tu+1}\right\|^2+\frac12\|v_x\|^2+\R \|\tu\|^2 \leq 0.
\end{align}
\vincent{Let us define}
    \begin{align}
        E(t)&:=\frac1\veps\intol (\eta(\tu(t)+1)-\eta(1)) +\frac 1{2\veps}\|v(t)\|^2+ \intol f(\tu+1)-f(1)\label{E:def}\\
        F(t)&=\frac12\|v_x(t)\|^2+\R \|\tu(t)\|^2.\label{G:def}
    \end{align}
\vincent{Then \eqref{extra} can be rewritten as}
\begin{align}
    \frac{d}{dt}E(t)+F(t) \leq 0.\notag
\end{align}
Let \vincent{$T^*$ denote the time asserted in \cref{bounds}. Using the fact that $\tu=u-1$, notice that we can rewrite \cref{bounds} as}
\begin{align}\notag
   \frac29 \|\tu(t)\|^2 &\leq \intol f(\tu(t)+1)-f(1) \leq 2\| \tu(t)\|^2\notag\\
     \frac13 \|\tu(t)\|^2 &\leq \intol \eta(\tu(t)+1)-\eta(1) \leq  \|\tu(t)\|^2,\notag
\end{align}
\vincent{for all $t\geq T^*$.} Using \vincent{these} bounds we get
\begin{align}\notag
    E(t)&\leq \frac1{2\veps} \|v(t)\|^2 +\left(2+\frac1\veps\right)\|\tu(t)\|^2\\
    &\leq\vincent{\frac{4}\veps} \left(\|\tu(t)\|^2+\|v(t)\|^2\right),\quad t\geq T^*.\notag
\end{align}
\vincent{Also, by the} Poincar\'e inequality, we have
\begin{align}\notag
    F(t)&\geq \vincent{\frac{r\wedge 1}2}\left(\|\tu(t)\|^2+\| v(t) \|^2\right),\quad t\geq T^*,
\end{align}
\vincent{where $a\wedge b:=\min\{a,b\}$. Hence}
\begin{align}\notag
    E(t)&\leq 
    \vincent{\frac{8}{ (r\wedge1)}}\veps F(t),\quad t\geq T^*.
\end{align}
Therefore, since $\veps\leq1$, we have
\begin{align}\notag
   \frac{d}{dt} E(t)+ \vincent{\frac{(r\wedge1)}8}\veps E(t)\leq 0,
\end{align}
\vincent{so that Gronwall's inequality gives}
\begin{align}\notag
   E(t)\leq E(T^*)e^{-\vincent{\frac{(r\wedge1)}8}\veps(t-T^*)},\quad t\geq T^*.
\end{align}
\vincent{Now from Lemma \eqref{bounds} we have that} for all $t\geq T^*$
\begin{align}
    E(t)&\geq \left(\frac1{3\veps} +\frac29\right)\|\tu(t)\|^2 +\frac1{2\veps}\|v(t)\|^2\notag \\
    &\geq \frac{1}{3\veps}\left(\|\tu(t)\|^2+\|v(t)\|^2\right).\notag
\end{align}
Therefore $\text{for all } t\geq T^*$
\begin{align}\label{est:L2:decay}
    \|\tu(t)\|^2+\|v(t)\|^2&\leq \vincent{3E(T^*)e^{-\frac{\veps r}8(t-T^*)}}.
\end{align}
\vincent{This establishes the asymptotic decay in $L^2$. Now we establish asymptotic decay in $H^1$}.

From \eqref{H1} and \cref{cor:H2:complete}, we have 
\begin{align}\label{uno'}
 \frac{d}{dt}&\left(\|\tu_x\|^2+\|v_x\|^2\right)+\|\tu_{xx}\|^2+\veps \|v_{xx}\|^2\notag \\
\leq & c(\|v_x\|^2+\|v\|^2+r^2\|\tu\|^2+\|\tu_x\|^2+\veps\|v_x\|^2)(\|\tu_x\|^2+\|v_x\|^2)+C(r,\veps)(\|\tu\|^2+\|v_x\|^2).\notag\\
\leq & C_1(r,\veps)\left(\|\tu\|^2+ \|\tu_x\|^2+\|v_x\|^2\right).
\end{align}
\vincent{Recall that from \eqref{extra} we have}
\begin{align}\label{extra2}
    \frac{d}{dt}\left(\frac1\veps\intol (\eta(\tu+1)-\eta(1)) +\frac 1{2\veps}\|v\|^2+ \intol (f(\tu+1)-f(1)) \right)
    +\frac1\veps\intol \frac{(\tu_x)^2}{\tu+1}+\left\|\frac{\tu_x}{\tu+1}\right\|^2+\frac12\|v_x\|^2+\R \|\tu\|^2 \leq 0.
\end{align}
\vincent{Using \eqref{cond:u:Linfty:squared}, we have $\frac23\leq(\tu(t) +1)^{-1}\leq 2$,  $\text{for all } t\geq T^*$, so that \eqref{extra2} becomes}
\begin{align}\label{dos}
    \frac{d}{dt}\left(\frac1\veps\intol (\eta(\tu+1)-\eta(1)) +\frac 1{2\veps}\|v\|^2+ \intol (f(\tu+1)-f(1)) \right)    +\vincent{\frac{2}{3\veps}\|\tu_x\|^2}+\frac49\left\|\tu_x\right\|^2+\frac12\|v_x\|^2+\R \|\tu\|^2 \leq 0.
\end{align}
{{}Define}
    \begin{align}
   {{} G(t)}=\|\tu_x\|^2+\|v_x\|^2+\frac{6C_1(r,\veps)}{(r\wedge1)\veps}\intol (\eta(\tu+1)-\eta(1)) +\frac{3C_1(r,\veps)}{(r\wedge1)\veps}\|v\|^2+ \frac{6C_1(r,\veps)}{(r\wedge1)}\intol (f(\tu+1)-f(1)).\notag
    \end{align}
Observe that by the {{}Poincar\'e inequality} and \cref{bounds}, we have
\begin{align}
    \|\tu_x\|^2+\|v_x\|^2\leq G(t)\leq {{}C_2(r,\veps)}(\|\tu_x\|^2+\|v_x\|^2+\|\tu\|^2),\label{bound:A}
\end{align}
{{}where $C_2(r,\veps)=c{(r\wedge 1)^{-1}\veps^{-1}}$}. \vincent{Now, by taking the product of $6C_1(1\wedge r)^{-1}$} with \eqref{dos} and adding \vincent{the result to} \eqref{uno'}, we get 
\begin{align}
    \frac{d}{dt}&G(t) +\veps\|v_{xx}\|^2+\|\tu_{xx}\|^2+C_1(\|\tu_x\|^2+\|v_x\|^2+\|\tu\|^2)\leq0.\label{tres}
\end{align}
{{}so that \eqref{bound:A} implies}
\begin{align}\notag
    \frac{d}{dt}&G(t)+{{}\frac{C_1}{C_2}}G(t)\leq0.
\end{align}
\vincent{An application of Gr\"onwall's inequality and the fact that $\eta,f$ are strongly convex yields}
\begin{align}\label{est:H1:decay}
   \vincent{ \|\tu_x(t)\|^2+\|v_x(t)\|^2 \leq e^{- {{}\frac{C_1(r,\veps)}{C_2(r,\veps)}}(t-T^*)}A(T^*)
   , \quad t\geq T^*,}
\end{align}
\vincent{which, coupled with \eqref{est:L2:decay}, establishes the asymptotic decay in $H^1$}. 


\vincent{Finally}, for the $H^2$ decay recall from \eqref{est:H2:ineq} and \cref{cor:H2:complete}, we have
\begin{align}
   \frac{d}{dt} (\|\tu_{xx}\|^2+\|v_{xx}\|^2)+&\|\tu_{xxx}\|^2+\veps\|v_{xxx}\|^2
    \leq C_3(\|\tu_x\|^2+\|\tu_{xx}\|^2+\|v_x\|^2+\|v_{xx}\|^2).\label{cuatro}
\end{align}
\vincent{Recall that $0<\veps<1$. We may also assume that $C_3\geq C_1\geq 1$. Upon multiplying \eqref{tres} by $2C_3(r,\veps)\veps^{-1}$, then adding the result to \eqref{cuatro}, we obtain}
\begin{align}
    \frac{d}{dt} \left(\frac{2C_3}{\veps}G(t)+\|\tu_{xx}\|^2+\|v_{xx}\|^2\right) + \|\tu_{xx}\|^2+\|v_{xx}\|^2+{{}\frac{C_3C_1}{\veps}}(\|\tu\|^2+\|\tu_x\|^2+\|v_x\|^2)\leq 0,\notag
\end{align}
which, by the assumptions on the constants and the inequality \eqref{bound:A}, we obtain
\begin{align}
    \frac{d}{dt} \left(\frac{2C_3}{\veps}G(t)+\|\tu_{xx}\|^2+\|v_{xx}\|^2\right) + {{}\frac{1}2\left(1\wedge\frac{C_1}{C_2}\right)}\left(\|\tu_{xx}\|^2+\|v_{xx}\|^2+{{}\frac{2C_3}{\veps}G(t)}\right)\leq 0.\notag
\end{align}
\vincent{A final application of Gr\"onwall's inequality yields}
\begin{align}
\vincent{\|\tu_{xx}(t)\|^2+\|v_{xx}(t)\|^2\leq {{}e^{-\frac{1}2(1\wedge \frac{C_1}{C_2})(t-T^*)}}\left(\frac{2C_3}{\veps}G(T^*)+\|\tu_{xx}(T^*)\|^2+\|v_{xx}(T^*)\|^2\right)}, \qquad t>T^*.\notag
\end{align}
\vincent{Combining this estimate with \eqref{est:L2:decay} and \eqref{est:H1:decay} completes the proof.}
\end{proof}
\section{Vanishing chemical diffusion limit: Proof of \cref{thm:diff}}\label{sect:diff}

\vincent{To prove \cref{thm:diff}, we will first establish estimates for strong solutions of \eqref{eq:main11}, \eqref{eq:main22} that are uniform in $\veps>0$. To do so, we will again develop a bootstrap from $L^2$ to $H^2$, but ensure at each step that the estimates are independent of $\veps>0$.}

\begin{nlem}[$\veps$--independent bounds]\label{exist.diff}
\vincent{Let $u_0,v_0\in H^2$ such that $\frac{du_0}{dx}, v_0\in H^1_0$ and $u_0\geq0$ and $M_0<\infty$. Let $(u,v)$ denote the unique, global strong solution of \eqref{eq:main11} corresponding to $\veps\in(0,1)$ and initial data $(u_0,v_0)$.} 
Then for all $T>0$
 \[
 \sup_{t\in[0,T]}\left(\|u(t)\|^2_{H^2}+\|v(t)\|^2_{H^2}\right)+\int_0^T\left(\|u(s)\|^2_{H^3}+\veps\|v(s)\|^2_{H^3}\right)ds\leq C(T)
 \]
 where $C$ depends only on $r,T$, but is independent of $\veps$.
\end{nlem}

\begin{proof}[Proof of \cref{exist.diff}]
Fix $T>0$. We initiate the bootstrap by obtaining $L^2$ estimates.

\flushleft{\textit{Step 1: $L^2$ estimates.}}  \vincent{From \eqref{L2} and the fact that $\tu+1\geq 0$, we may apply Gronwall's inequality to obtain}
\begin{align}\label{est:L2:diff}
    \sup_{t\in[0,T]}\left(\|\tu(t)\|^2+\|v(t)\|^2 +\intot \|\tu_x(s)\|^2ds+\veps\intot\|v_x(s)\|^2ds \right)\leq e^{cT}\|\tu_0\|^2\leq C(T),
\end{align}
\vincent{for some constant $C$ that depends on $T$, but is independent of $\veps$.}

\flushleft{\textit{Step 2: $H^1$ estimates.}}
\vincent{We apply $\partial_x$ to \eqref{shift:eqn:u} and \eqref{shift:eqn:v}, then take the $L^2$ inner product with $\tu_x$ and $v_x$, respectively, and sum the results. After integrating by parts and applying \cref{rk:extrabc} we get}
\begin{align}\label{eq:H1:balance:diff}
    \frac 12\frac d{dt}(\|\tu_x\|^2+\|v_x\|^2)+\|\tu_{xx}\|^2+\veps\|v_{xx}\|^2\notag
    &= \intol (\tu v)_x\tu_{xx} +r\intol\tu_x(\tu(1+\tu))_x-\veps\intol(v^2)_{x}v_{xx}\notag\\
    &=\intol (\tu v)_x\tu_{xx}+r\intol \tu_{x}^2(1+2\tu)-2\veps\intol vv_xv_{xx}\notag\\
    &=I_1+I_2+I_3.
\end{align}
\vincent{By H\"older's inequality, \eqref{GNS}, \eqref{est:L2:diff}, and Young's inequality, we get}
\begin{align}
    |I_1| &\leq \frac 12 \|\tu_{xx}\|^2 +\|\tu\|_\infty^2\|v_x\|^2+\|v\|_\infty^2\|\tu_x\|^2\notag\\
    &\leq \frac 12\|\tu_{xx}\|^2 +c(\|\tu\|^2+\|\tu_x\|^2)\|v_x\|^2 +c(\|v_x\|^2+\|v\|^2)\|\tu_x\|^2\notag\\
    &\leq \frac 12\|\tu_{xx}\|^2 +c\|\tu_x\|^2\|v_x\|^2+C(T) (\|\tu_x\|^2+\|v_x\|^2).\notag
\end{align}
\vincent{Similarly, we estimate}
\begin{align}
    |I_2| &\leq  r\|1+2\tu\|_\infty \|\tu_x\|^2 \notag\\
    &\leq r\|\tu_x\|^2 +c\|\tu\|_\infty\|\tu_x\|^2\notag\\
    &\leq r\|\tu_x\|^2 +c(\|\tu\|^{1/2}\|\tu_x\|^{1/2}+\|\tu\|)\|\|\tu_x\|^2\notag\\
    &\leq C(T)\|\tu_x\|^2 +c(\|\tu\|^{2/3}+\|\tu_x\|^{2})\|\|\tu_x\|^2\notag\\
    &\leq C(T)\|\tu_x\|^2 +c\|\tu_x\|^{2}\|\tu_x\|^2.\notag
\end{align}
\vincent{Finally, additionally using the fact that $\veps<1$, we get}
\begin{align}
    |I_3|&\leq \frac\veps 2\|v_{xx}\|^2 +c\veps\|v\|_\infty^2\|v_x\|^2\notag\\
    &\leq \frac\veps 2\|v_{xx}\|^2 +c\veps(\|v\|^2+\|v_x\|^2)\|v_x\|^2\notag\\
    &\leq \frac\veps 2\|v_{xx}\|^2 +C(T)\veps\|v_x\|^2+c\veps\|v_x\|^2\|v_x\|^2\notag\\
    &\leq \frac\veps 2\|v_{xx}\|^2 +C(T)\|v_x\|^2+c\veps\|v_x\|^2\|v_x\|^2.\notag
\end{align}
\vincent{Combining $I_1$--$I_3$ in \eqref{eq:H1:balance:diff}, we get} 
\begin{align}
  \vincent{\frac d{dt}}(\|\tu_x\|^2+\|v_x\|^2)+\|\tu_{xx}\|^2+\veps\|v_{xx}\|^2
  \leq c(\|\tu_x\|^2+\veps\|v_x\|^2)(\|\tu_x\|^2+\|v_x\|^2)+C(T)(\|\tu_x\|^2+\|v_x\|^2).\label{est:H1:balance:diff}
\end{align}
\vincent{An application of Gr\"onwall's inequality and \eqref{est:L2:diff}, yields}
\begin{align}\label{est:H1:diff}
  \sup_{t\in[0,T]}\left(\|\tu_x(t)\|^2+\|v_x(t)\|^2+\intot\|\tu_{xx}(s)\|^2+\veps\intot\|v_{xx}(s)\|^2\right)\leq C(T).
\end{align}

\textit{Step 3: $H^2$ estimates.} \vincent{We apply $\partial_{xx}$ to \eqref{shift:eqn:u} and \eqref{shift:eqn:v}, then take the $L^2$ inner product with $\tu_{xx}$ and $v_{xx}$, respectively, and sum the results. After integrating by parts and applying \cref{rk:extrabc} we get}
\begin{align}
\frac 12\frac{d}{dt}\left( \|\tu_{xx}\|^2+\|v_{xx}\|^2\right)&+\|\tu_{xxx}\|^2+\veps\|v_{xxx}\|^2\notag\\
&=\intol (\tu v)_{xx}\tu_{xxx}+\R \intol(2\tu+1)\tu_{x}\tu_{xxx}-\veps\intol(v^2)_{xx}v_{xxx}\notag\\
&=\intol (\tu_{xx}v+2\tu_x v_x+\tu v_{xx})\tu_{xxx}+\R \intol(2\tu+1)\tu_{x}\tu_{xxx}-\veps\intol(v^2)_{xx}v_{xxx}\notag\\
&= J_1+J_2+J_3.\label{eq:H2:balance:diff}
\end{align}
\vincent{To estimate the terms in the right hand side, we use H\"older's inequality, \eqref{GNS}, \eqref{est:L2:diff}, \eqref{est:H1:diff}, and Young's inequality. For $J_1$, we estimate}
\begin{align}
    |J_1|
    &\leq \frac14\|\tu_{xxx}\|^2 +\|v\|^2_\infty\|\tu_{xx}\|^2 +2\|\tu_x\|^2_\infty\|v_x\|^2 +\|\tu\|^2_\infty\|v_{xx}\|^2\notag\\
    &\leq \frac14\|\tu_{xxx}\|^2 +c(\|v\|^2+\|v_{xx}\|^2)\|\tu_{xx}\|^2 +c(\|\tu_x\|^2+\|\tu_{xx}\|^2)\|v_x\|^2 +c(\|\tu\|^2+\|\tu_x\|^2)\|v_{xx}\|^2\notag\\
    &\leq \frac14\|\tu_{xxx}\|^2 +C(T)(\|\tu_{xx}\|^2 +\|v_{xx}\|^2)+C(T).\notag
\end{align}
Similarly, we estimate
\begin{align}
    |J_2|&\leq \frac14\|\tu_{xxx}\|^2+r^2\|1+2\tu\|^2_\infty\|\tu_x\|^2\notag\\
    &\leq\frac14\|\tu_{xxx}\|^2+C(T)\|\tu_x\|^2\notag\\
    &\leq\frac14\|\tu_{xxx}\|^2+C(T).\notag
\end{align}
 \vincent{Lastly, additionally using the fact that $\veps<1$, we estimate}
\begin{align}
    |J_3|&\leq \frac\veps2\|v_{xxx}\|^2 +2(\veps\|v\|_\infty^2\|v_{xx}\|^2+\veps\|v_x\|^2\|v_x\|^2)\notag\\
    &\leq  \frac\veps2\|v_{xxx}\|^2 +c\veps\|v\|^2\|v_{xx}\|^2+c\veps\|v_x\|^2\|v_{xx}\|^2+c\veps\|v_x\|^2\|v_x\|^2\notag\\
    &\leq  \frac\veps2\|v_{xxx}\|^2 +C(T)\|v_{xx}\|^2+c\|v_x\|^2\|v_{xx}\|^2+c\|v_x\|^2\|v_x\|^2\notag\\
    &\leq  \frac\veps2\|v_{xxx}\|^2 +C(T)\|v_{xx}\|^2+C(T).\notag
\end{align}
\vincent{Combining $J_1$--$J_3$ in \eqref{eq:H2:balance:diff}}, we get
\begin{align}\label{est:H2:balance:diff}
  \frac{d}{dt}\left( \|\tu_{xx}\|^2+\|v_{xx}\|^2\right)+\|\tu_{xxx}\|^2+\veps\|v_{xxx}\|^2\leq   C(T)(\|\tu_{xx}\|^2 +\|v_{xx}\|^2)+C(T),
\end{align}
\vincent{so that Gronwall's inequality implies}
\begin{align}\label{est:H2:diff}
  \left( \|\tu_{xx}\|^2+\|v_{xx}\|^2\right)&+\intot(\|\tu_{xxx}\|^2+\veps\|v_{xxx}\|^2)\leq   C(T),
\end{align}
as desired.
\end{proof}


Now we are ready to prove \cref{thm:diff}. 
\begin{proof}[Proof of \cref{thm:diff}]
\vincent{Let $\veps\geq0$ and $(\tu^\veps,v^\veps)$ denote the unique strong solution to \eqref{shift:eqn:u}, \eqref{shift:eqn:v}, \eqref{ic1}, \eqref{NBC1}. Let}
\begin{align}\notag
   {U}=\tu^\veps-\tu^0,\qquad 
   {V}=v^\veps-v^0.
\end{align}
We obtain the following system for \vincent{$({U},{V})$:}
\begin{align}
{U}_t+(\tu^0 {V}+{U} v^\veps)_x+{V}_x&={U}_{xx}-r{U}({U}+1+2\tu^0)\label{eqn1e}\\
{V}_t+{U}_x&=\varepsilon  v^\veps_{xx}+\varepsilon((v^\veps)^2)_x\label{eqn2e}
\end{align}
with initial conditions given by
\begin{align}\notag
{U}_0=0={V}_0,\quad ({U}_0,{V}_0)\in H^2(\Omega)
\end{align}
and boundary conditions given by
\begin{align}\notag
{U}_x|_{\partial\Omega}=0, \quad v^\veps|_{\partial\Omega}=0.
\end{align}

\vincent{Upon taking the $L^2$ inner product of \eqref{eqn1e}, \eqref{eqn2e} with ${U}, {V}$, respectively, summing the results, then using integration by parts and the boundary conditions}, we get 
\begin{align}\label{eq:est:L2:diff:lim}
\frac 12 \frac{d}{dt}& \left( \|{U}\|^2 +\| {V} \|^2 \right) +\|{U}_x\|^2\notag\\
=&-\intol  ( \tu^0 {V}+{U} v^\veps)_x {U} -r \intol {U}^2({U}+1+2\tu^0)+\veps\intol v_{xx}^\veps {V} +\veps \intol \bdy_x(v^\veps)^2{V} \notag \\
=&-\intol  (\tu^0_x {V}+\tu^0{V}_x-{U}_xv^\veps){U}  -r \intol ({U}^3-{U}^2 -2 {U}^2\tu^0)+\veps\intol v_{xx}^\veps {V}+2\veps\intol v^\veps v^\veps_x{V}\notag\\
=&K_1+K_2+K_3+K_4.
\end{align}
Using \eqref{GNS}, and \eqref{est:L2:diff} and \eqref{est:H1:diff}, we estimate $K_1$ as
\begin{align}
|K_1|
&\leq  c\|\tu^0_x\|^2_\infty\|{U}\|^2+c\|{V}\|^2+c\|\tu^0\|_\infty^2\|{U}\|^2+ c\|{V}_x\|^2 +\frac 12 \|U_x\|^2+c\|v^\veps\|^2_\infty\|{U}\|^2)\notag\\
&\leq  \frac 12 \|{U}_x\|^2 +C(T)( \|{U}\|^2+\|{V}\|^2) +c\|{V}_x\|^2\notag.
\end{align}
\vincent{Making use of \eqref{GNS}, \eqref{est:L2:diff}, \eqref{est:H1:diff} and H\"older inequality, we estimate $K_2$ as}
\begin{align}
|K_2|
\leq &cr\|{U}\|_\infty\|{U}\|^2+cr \|{U}\|^2+c r\|\tu^0\|_\infty\|{U}\|^2 \notag\\
\leq &cr(\| \tu^\veps\|_\infty+\|\tu^0\|_\infty)\|{U}\|^2+cr \|{U}\|^2+cr \|\tu^0\|_\infty\|{U}\|^2 \notag\\
\leq& C(T)\|{U}\|^2.\notag
\end{align}
 \vincent{Similarly, we estimate $K_3$ as}
\begin{align}
|K_3|
\leq \frac{\veps^2}{2}\|v_{xx}^\veps\|^2+\frac 12 \|{V}\|^2\notag.
\end{align}
\vincent{Lastly, we estimate $K_4$ as}
\begin{align}
|K_4|
&\leq \veps^2\|v^\veps\|_\infty^2\|v^\veps_x\|^2+\|{V}\|^2\notag.
\end{align}
\vincent{Combining $K_1$--$K_4$ in \eqref{eq:est:L2:diff:lim}}, we have
\begin{align}\label{est:L2:diff:lim}
\frac{d}{dt}&\left( \|{U}\|^2 +\| {V}\|^2 \right) +\frac12\|{U}_x\|^2 \notag\\
\leq & C(T)((1\vee r)\|{U}\|^2+\|{V}\|^2) +C\|{V}_x\|^2+\frac{\veps^2}{2}\|v_{xx}^\veps\|^2+\veps^2\|v^\veps\|_\infty^2\|v^\veps_x\|^2.
\end{align}

\vincent{We take $\partial_x$ of \eqref{eqn1e}, \eqref{eqn2e}, then taking the $L^2$ inner product of ${U}_x,{V}_x$, respectively, and adding the results. After using integration by parts and the initial conditions we obtain}
\begin{align}
\frac12&\frac{d}{dt} (\|{U}_x\|^2 +\|{V}_x\|^2 ) +\|{U}_{xx}\|^2 \notag\\
=&\intol (\tu^0 {V}+{U} v^\veps)_{x}{U}_{xx}-r\intol {U}({U} +2\tu^0+1){U}_{xx}-\intol{U}_{xx}{V}_x+\veps\intol v^\veps_{xxx}{V}_x+2\veps\intol(v^\veps v^\veps_x)_x{V}_x\notag\\
=&\intol (\tu^0_x {V}+\tu^0 {V}_x+{U}_x v^\veps+{U} v^\veps_x){U}_{xx}-r\intol{U}({U} +2\tu^0+1){U}_{xx}-\intol{U}_{xx}{V}_x\notag\\
&+\veps\intol v^\veps_{xxx}{V}_x+2\veps\intol ((v^\veps_x)^2+v^\veps v^\veps_{xx}){V}_x\notag\\
=& L_1+L_2+L_3+L_4+L_5. \label{eq:H2:diff:lim}
\end{align}
\vincent{We proceed to estimate the right-hand side terms using H\"older's inequality, \eqref{GNS}, \eqref{est:L2:diff}, \eqref{est:H1:diff}, and Young's inequality. We estimate $L_1$ as}
\begin{align}
|L_1|
&\leq \frac16\|U_{xx}\|^2 +\|\tu^0_x\|^2_\infty\|{V}\|^2+\|\tu^0\|_\infty^2\|{V}_x\|^2+\|v^\veps\|^2_\infty\|{U}\|^2+\|v_x^\veps\|^2_\infty\|{U}_x\|^2\notag\\
&\leq \frac16\|U_{xx}\|^2 +C(T)(\|{V}\|^2+\|{V}_x\|^2+\|{U}\|^2+\|{U}_x\|^2).\notag
\end{align}
\vincent{We estimate $L_2$ as}
\begin{align}
|L_2|
&\leq \frac16\|{U}_{xx}\|^2+c\|{U}\|^2_\infty\|{U}\|^2+c\|{U}\|^2+c\|\tu^0\|^2_\infty\|{U}\|^2\notag\\
&\leq \frac16\|{U}_{xx}\|^2+C(T)\|{U}\|^2.\notag
\end{align}
\vincent{We estimate $L_3$ as}
\begin{align}
|L_3|&\leq \frac 16\|{U}_{xx}\|^2+c\|{V}_x\|^2.\notag
\end{align}
\vincent{We estimate $L_4$ as}
\begin{align}
|L_4|&\leq \frac{\veps^2}2\|v^\veps_{xxx}\|^2+\frac12\|{V}_x\|^2.\notag
\end{align}
\vincent{We estimate $L_5$ as}
\begin{align}
|L_5|
&\leq \frac12 \|{V}_x\|^2 +\frac{\veps^2}2\|v^\veps_x\|^2_\infty\|v^\veps_{x}\|^2+\frac{\veps^2}2\|v^\veps\|^2_\infty\|v_{xx}^\veps\|^2.\notag
\end{align}
\vincent{Combining $L_1$--$L_5$ in \eqref{eq:H2:diff:lim}}, we get
\begin{align}
\frac{d}{dt}&(\|{U}_x\|^2+\|{V}_x\|^2)+ \|{U}_{xx}\|^2\notag\\
&\leq C(T)(\|{V}\|^2+\|{V}_x\|^2+\|{U}\|^2+\|{U}_x\|^2 )+{\veps^2}\|v^\veps_{xxx}\|^2+{\veps^2}\|v^\veps_x\|^2_\infty\|v^\veps_{x}\|^2+{\veps^2}\|v^\veps\|^2_\infty\|v_{xx}^\veps\|^2.\label{est:H2:diff:lim1}
\end{align}
\vincent{After adding \eqref{est:H2:diff:lim1} to \eqref{est:L2:diff:lim}, we arrive at}
\begin{align}
\frac{d}{dt}&(\|{U}\|^2+\|{V}\|^2+\|{U}_x\|^2+\|{V}_x\|^2)+ \|{U}_{x}\|^2+ \|{U}_{xx}\|^2\notag\\
&\leq C(T)(\|{V}\|^2+\|{V}_x\|^2+\|{U}\|^2+\|{U}_x\|^2 )+C(T)\veps^2 (\|v_{x}^\veps\|^2+\|v_{xx}^\veps\|^2+\|v_{xxx}^\veps\|^2).\notag
\end{align}
From \eqref{est:L2:diff}, \eqref{est:H1:diff}, and \eqref{est:H2:diff} we know that 
\begin{align}
\veps^2\intot (\|v_{x}^\veps\|^2+\|v_{xx}^\veps\|^2)+\veps^2\intot\|v_{xxx}^\veps\|^2\leq  C(T)\veps\notag. 
\end{align}
Therefore, \vincent{after applying Gronwall's inequality}, we obtain
\begin{align}
\|{U}(t)\|^2+\|{V}(t)\|^2+\|{U}_x(t)\|^2+\|{V}_x(t)\|^2+\intot( \|{U}_{x}(s)\|^2+ \|{U}_{xx}(s)\|^2)ds
&\leq C(T)\veps,\notag
\end{align}
\vincent{for some constant $C(T)$, depending on $T$, but independent of $\veps$. As this holds for all $0\leq t\leq T$, upon passing the limit $\veps\rightarrow0$, we establish the claim.}
\end{proof}

\section{Numerical results}

\begin{figure}[htpb]
    \centering
\boxed{\includegraphics[scale=0.2]{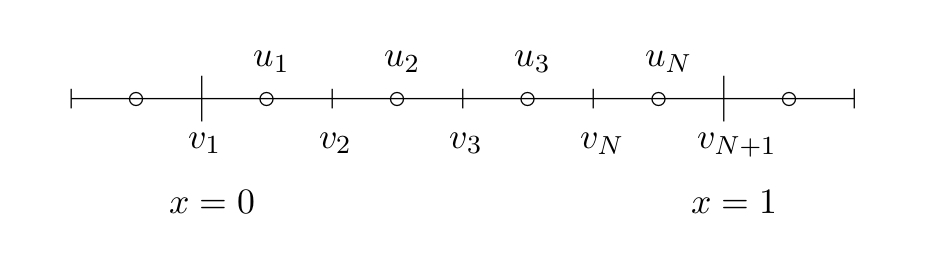}}
\caption{Diagram for the splitting of the domain and evaluation of $u$ and $v$}
\label{fig:domain}
\end{figure}

\begin{figure}[htbp]
    \centering
    \begin{subfigure}[h]{0.43\textwidth}
        \includegraphics[width=\textwidth]{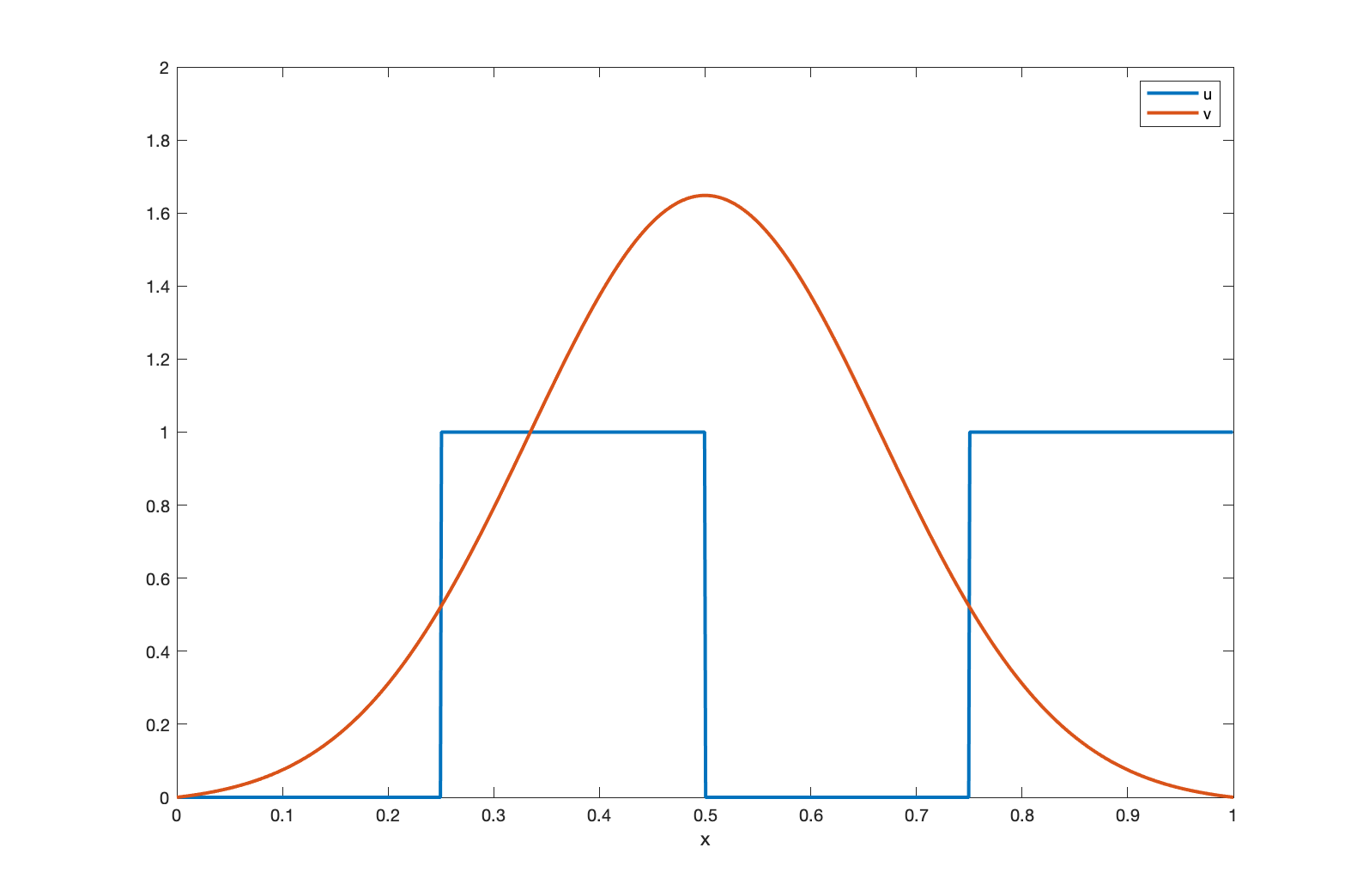}
        \caption{}
        \label{fig:initial}
    \end{subfigure}
    \qquad
    \begin{subfigure}[h]{0.43\textwidth}
        \includegraphics[width=\textwidth]{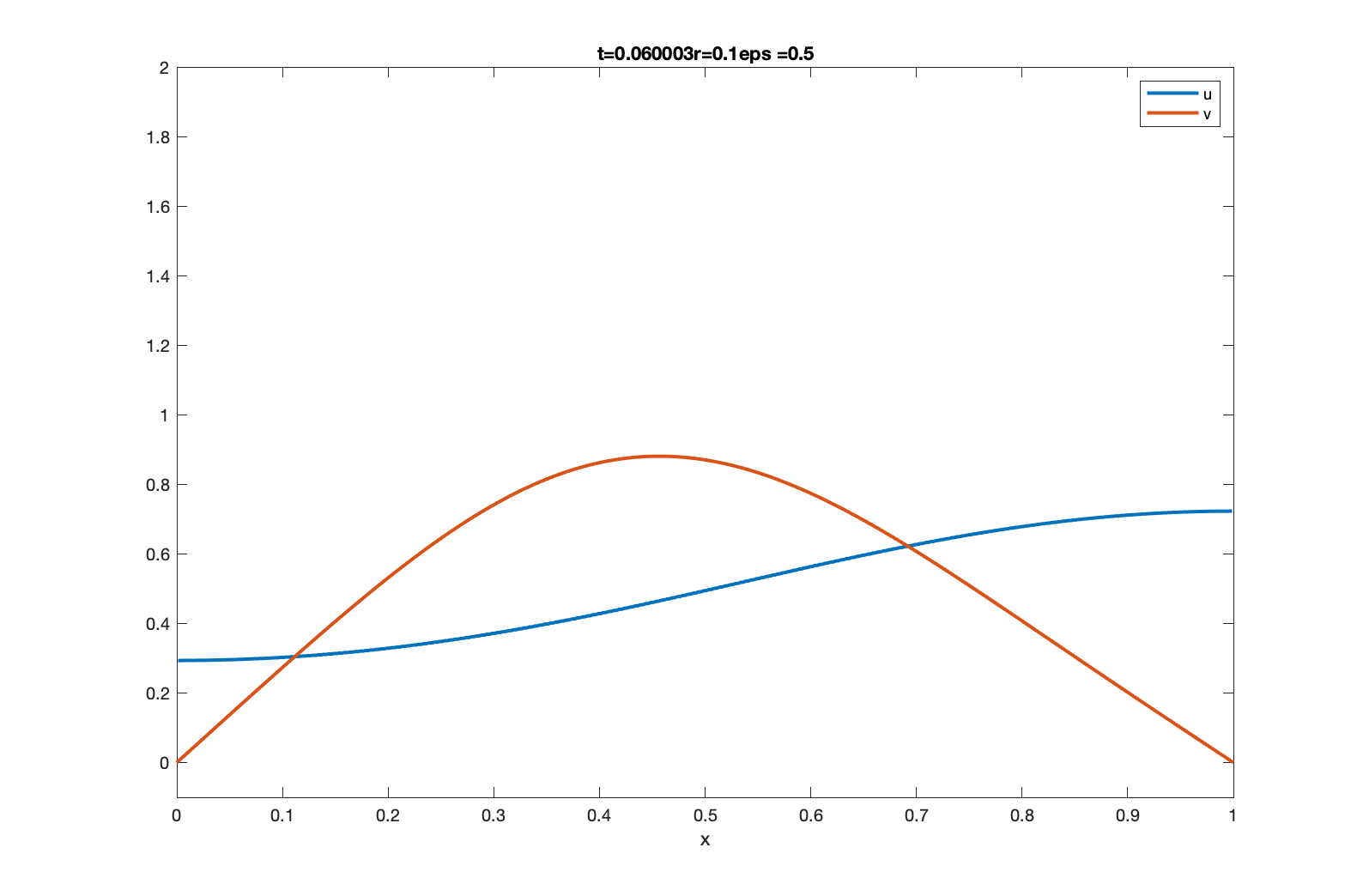}
        \caption{}
        \label{fig:transient:a}
    \end{subfigure}
    
    \begin{subfigure}[h]{0.43\textwidth}
        \includegraphics[width=\textwidth]{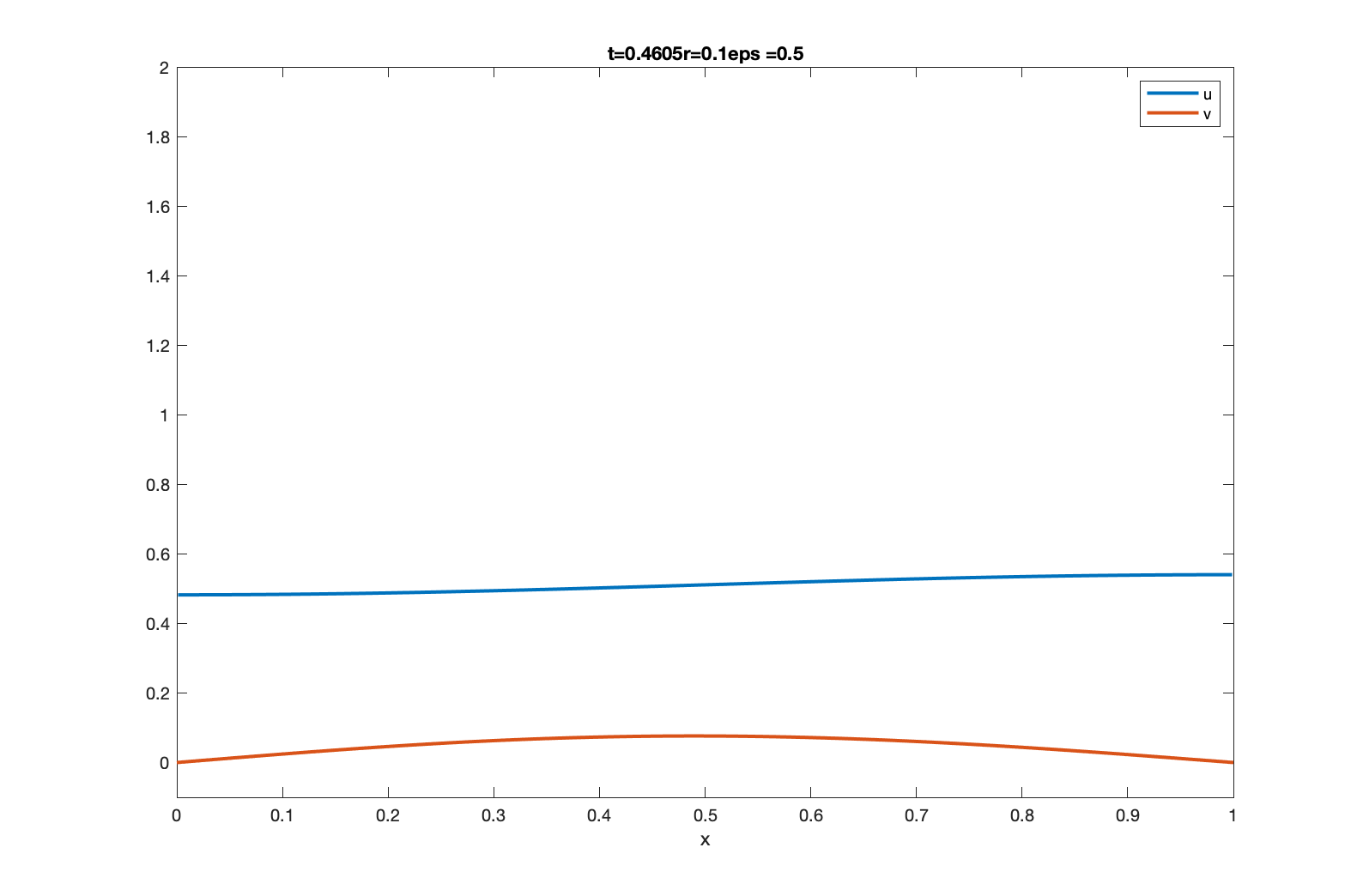}
        \caption{}
        \label{fig:transient:b}
    \end{subfigure}
    \qquad
    \begin{subfigure}[h]{0.43\textwidth}
        \includegraphics[width=\textwidth]{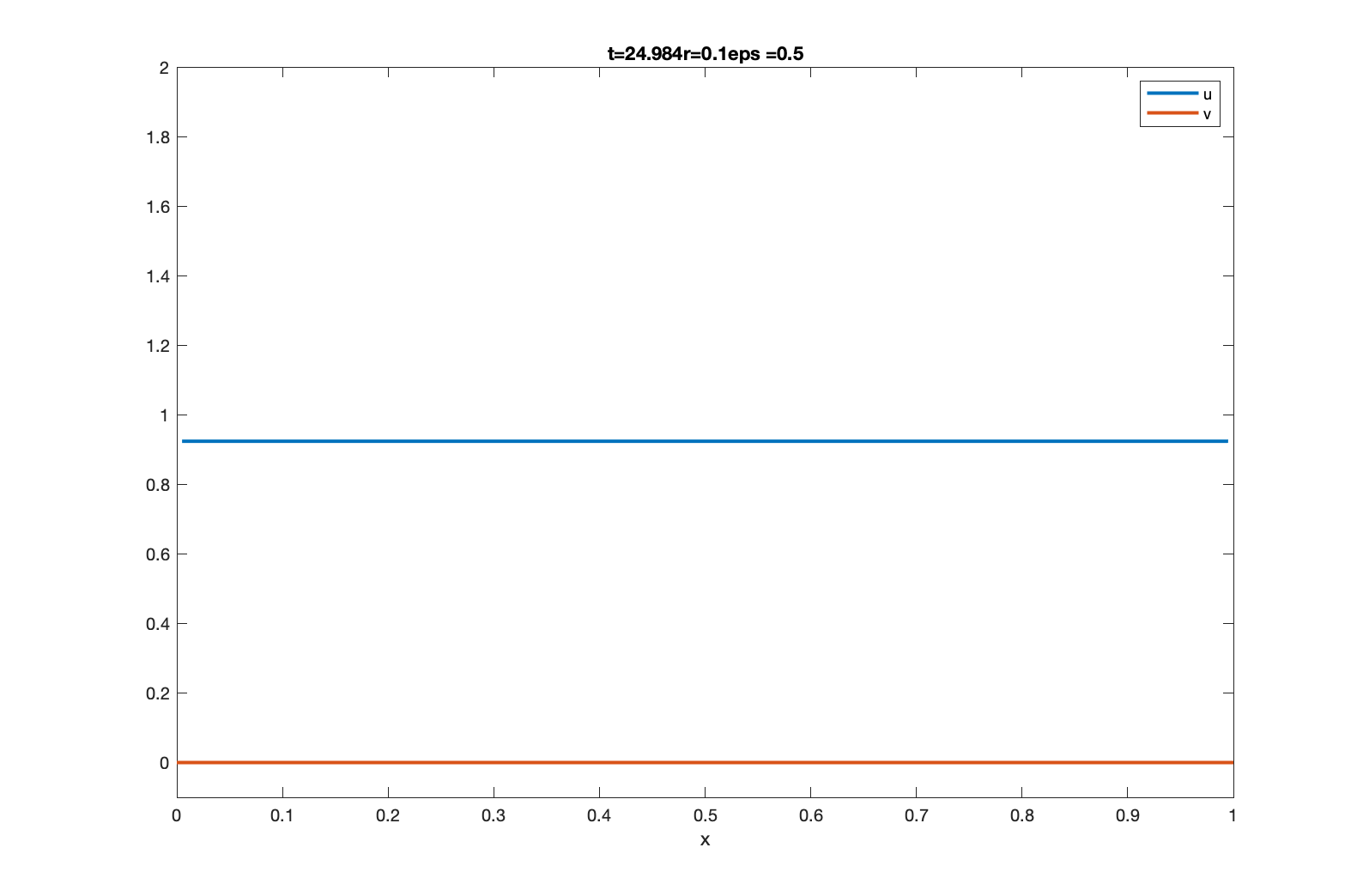}
        \caption{}
        \label{fig:logistic}
    \end{subfigure}
    \caption{\cref{fig:initial} shows the initial conditions for $u$ and $v$. The rest show the behaviour of the solutions for $\veps=0.5$ and $r=0.1$ at times $t=0.06$, $t=0.46$ and $t=24.98$} 
    \label{fig:behavior}
\end{figure}

\vincent{In this section, we carry out various numerical tests to study the dynamical properties of solutions to \eqref{eq:main11}, \eqref{eq:main22}. In particular, we 1) provide numerical confirmation of the rigorous qualitative results established above (cf. \cref{thm:Glob:exist} and \cref{thm:diff}), 2) demonstrate the phenomenon of separation of scales between the diffusive regime and logistic regime, and 3) identify a robust transient behavior of the solutions that was not treated by the mathematical analysis above.}
\vincent{The chemotaxis model \eqref{eq:keller:segel} with logarithmic sensitivity is generally difficult to solve using routine numerical methods due to the singularity of the term $\frac{\ln c}{c}$ appearing in the chemical concentration. For this reason, we instead provide simulations for the transformed system \eqref{eq:main11}, \eqref{eq:main22}. To carry out our numerical tests, we employ an explicit finite-difference scheme} to solve the equations with a second order approximation of \vincent{spatial} derivatives. For this reason, the \vincent{temporal mesh was chosen to prescribe by} $\Delta t= \frac{(\Delta x)^2}2$, similar to what would be chosen for the heat equation. \vincent{To ensure that numerical diffusion does not dominate chemical diffusion when $\veps\ll 1$, we use a spatial mesh size} of $\Delta x<\sqrt{\frac{\veps}{10}}$. The domain for the \vincent{system is normalized to  be the interval $[0,1]$ with two floating points to define the Neumann boundary conditions on $u$ by imposing $u(-1)=u(1)$ and $u_x(0)=u(1)-u(-1)$ (see \cref{fig:domain}). For the following results, we have used initial data $u(x,0)=\text{Heaviside}(x-0.25)-\text{Heaviside}(x-0.5)+\text{Heaviside}(x-0.75)$ and $v(x,0)=\frac 53e^{-18(x-\frac12)^2}$ (see \cref{fig:initial}). Note that the spatial average of $u_0$ is $0.5$. This particular choice of initial conditions was made because it exhibited non-trivial behavior in comparison with sinusoidal initial data, which appeared to rapidly relax to a constant. We also point out that in our simulations, the initial datum satisfies $u_0=0$ on two disjoint intervals, so that our choice of initial data actually lies outside of the conditions specified in our rigorous theorems above (see \cref{thm:Lyapunov}, \cref{thm:Glob:exist}, \cref{thm:diff}). Indeed, there we assumed that $u_0$ satisfies particular integral conditions on its logarithm (see \eqref{def:M}). This assumption is very natural for the analysis and is commonly assumed in the literature. Nevertheless, we do not observe any significant change of behavior by allowing for the initial data to be identically zero on an interval, thus suggesting that the assumption made in \eqref{def:M} is purely a technical one for the analysis. It would be interesting to study if one can remove this assumption from the results above in a future work.}

\newpage
\begin{figure}[htbp]
\includegraphics[scale=0.4]{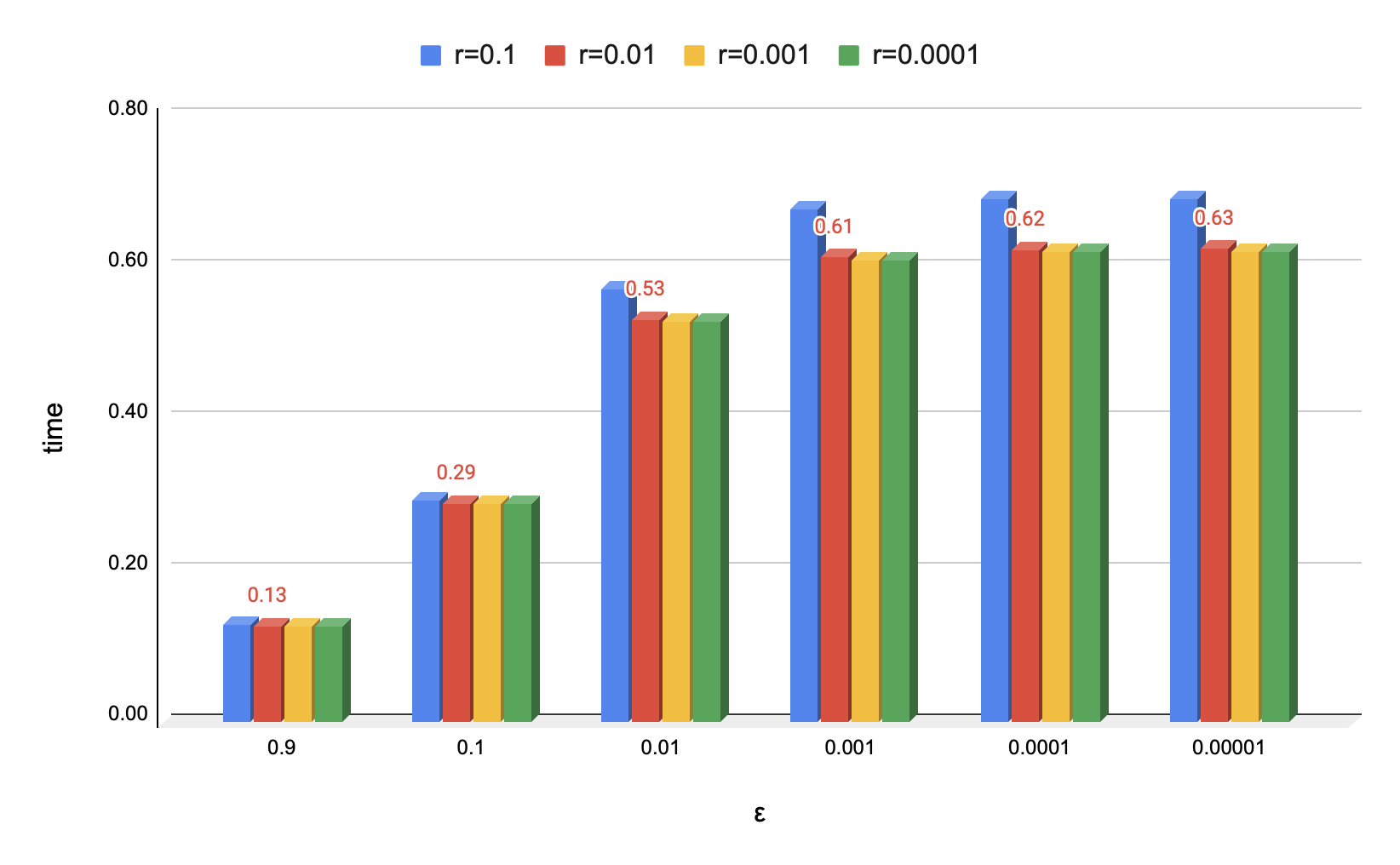}
        \caption{Time of relaxation to the initial average}
        \label{fig:ini:ave}
\end{figure}   
\begin{figure}[htbp]
\includegraphics[scale=0.4]{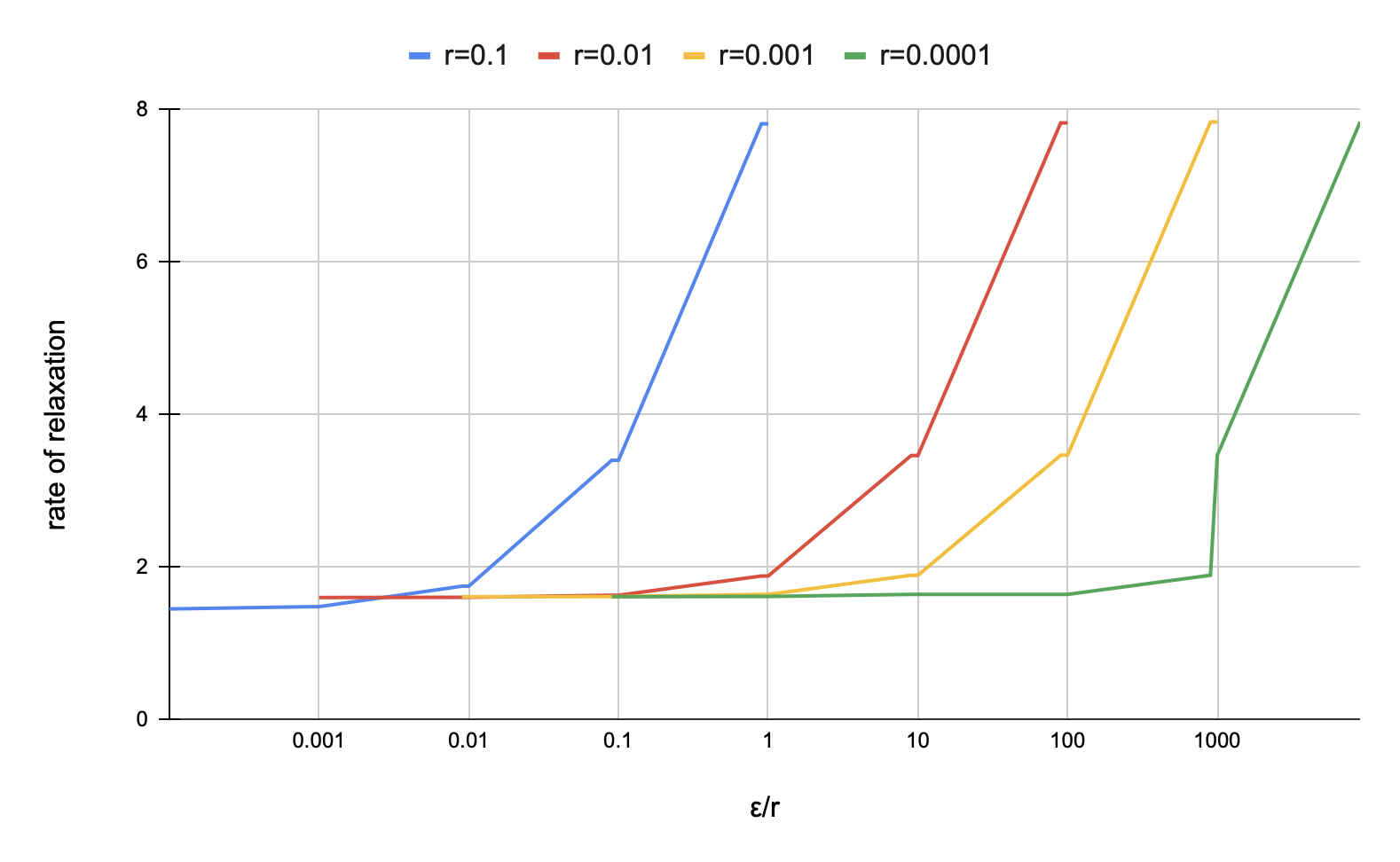}
        \caption{Relaxation rate to the initial average}
         \label{fig:relax:ave}
\end{figure} 

\vincent{In \cite{li2012global}, Li, Pan and Zhao showed that in the case of non-logistic growth, $r=0$, $u$ asymptotically converges in time to its initial average. In contrast to this behavior, in \cref{thm:Glob:exist} we proved that when $r>0$, then $u$ instead converges to its carrying capacity in the time-asymptotic limit, thus suggesting the eventual domination of the logistic mechanism over diffusive or chemotactic mechanisms. In numerically verifying this behavior, we also observe a robust transient behavior of the solution, in which $u$ \textit{first} tends to approach its initial average, before transitioning to the asymptotic behavior described by \cref{thm:Glob:exist} that characterizes the logistic dynamics (see \cref{fig:behavior}). Indeed, in \cref{fig:transient:a}, we first observe the initial data relaxing to the initial average before it begins to transition out of this behavior in \cref{fig:transient:b}, and march towards its final destination, i.e., the carrying capacity, in \cref{fig:logistic}. Therefore, we observe a separation of scales phenomenon for this model. In particular, the simulations indicate that the diffusive time scale is initially dominant and the solution exhibits non-logistic behavior. After a transient time, however, the logistic time scale becomes dominant and the dynamics change accordingly.}

\begin{figure}[htbp]
        \includegraphics[scale=0.5]{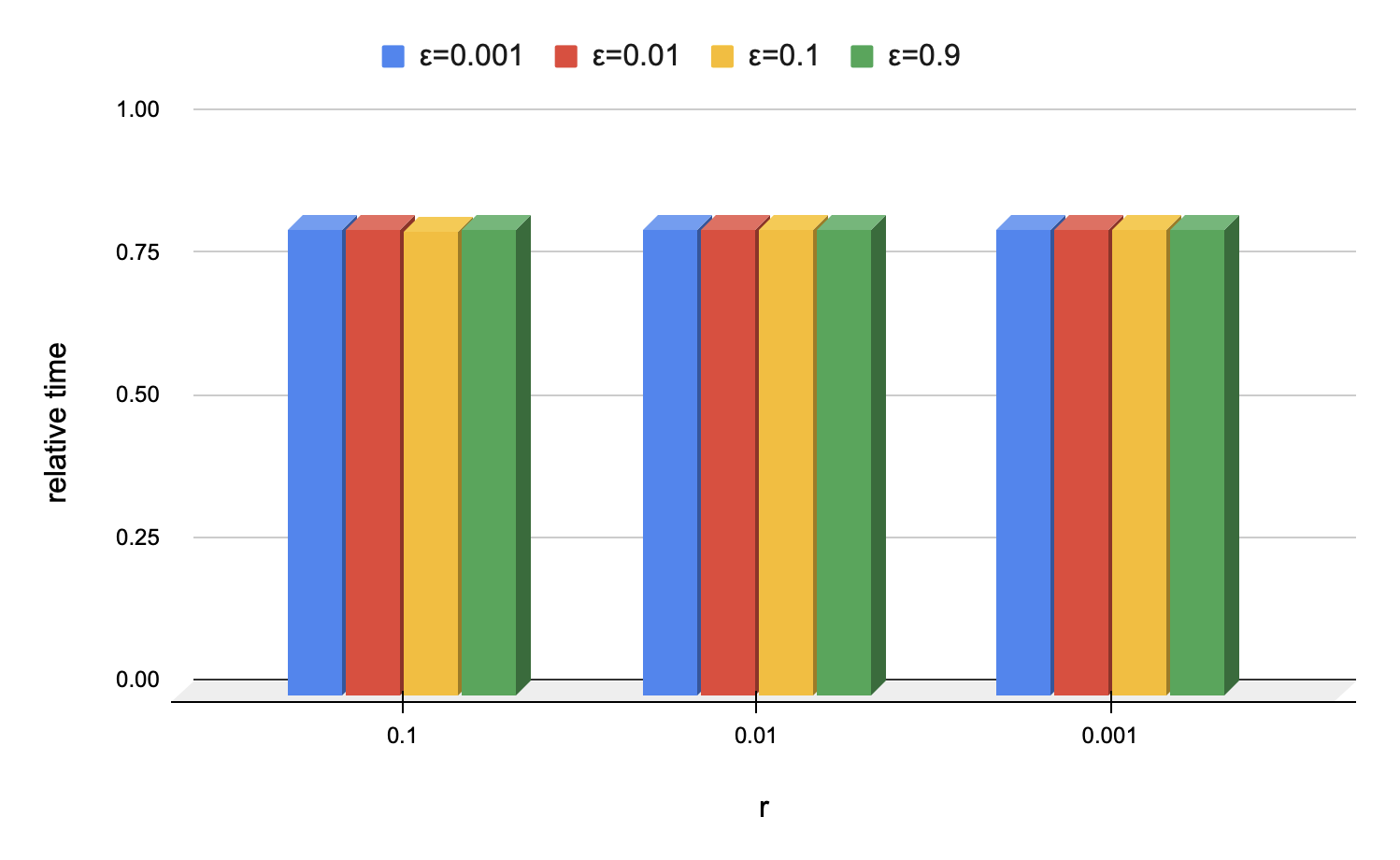}
        \caption{Relative time to leave the initial average and getting to the carrying capacity}
        \label{fig:rel:time}
 \end{figure} 

\begin {table}[htbp]
\centering
\small{\begin{tabular}{|c|c|c|c|c|c|c|c|}
 \hline
 \multicolumn{4}{|c|}{\small{\cellcolor{palegolden} \textbf{Time to leave initial average}}}&\multicolumn{4}{|c|}{\cellcolor{palesilver}\small{\textbf{Time to get to carrying capacity}}} \\
 \hline
 $\varepsilon$ & $r=1$ &\multicolumn{2}{|c|}{ $r=0.1$}&\multicolumn{2}{|c|}{ $r=0.01$}&\multicolumn{2}{|c|}{$r=0.001$}\\
 \hline
 $0.0001$ & \cellcolor{palesilver}2.4187    &\cellcolor{palegolden} 4.0861 &  \cellcolor{palesilver}22.058 & \cellcolor{palegolden}40.6292 & \cellcolor{palesilver}219.8052 & \cellcolor{palegolden}
405.55 & \cellcolor{palesilver}2197.3\\
 \hline
 $  0.001$ &\cellcolor{palesilver}2.4186   & \cellcolor{palegolden}4.0893 & \cellcolor{palesilver} 22.0566 & \cellcolor{palegolden}40.6279 & \cellcolor{palesilver}219.8038 & \cellcolor{palegolden}405.55 &\cellcolor{palesilver} 2197.3\\
 \hline
$ 0.01$  & \cellcolor{palesilver}2.3636   &\cellcolor{palegolden} 4.1058 &  \cellcolor{palesilver}22.0457 & \cellcolor{palegolden}40.6176 & \cellcolor{palesilver}219.7936 & \cellcolor{palegolden}405.53 & \cellcolor{palesilver}2197.3\\
 \hline
 $ 0.1$   & \cellcolor{palesilver}2.2427    & \cellcolor{palegolden}4.0927 & \cellcolor{palesilver} 22.0103 & \cellcolor{palegolden}40.58 & \cellcolor{palesilver}219.760 & \cellcolor{palegolden}405.50 & \cellcolor{palesilver}2197.3\\
 \hline
 $ 0.9$ & \cellcolor{palesilver}2.2145   & \cellcolor{palegolden}4.07 &  \cellcolor{palesilver}21.9894 & \cellcolor{palegolden}40.56 & \cellcolor{palesilver}219.7396 & \cellcolor{palegolden}405.48 & \cellcolor{palesilver}2197.3\\
 
 \hline
\end{tabular}}
\caption{\small{Time for $u$ to leave the initial average and to get to its carrying capacity}} \label{fig:table}

\end{table}

\begin{figure}[htbp]
   \qquad \includegraphics[scale=0.45]{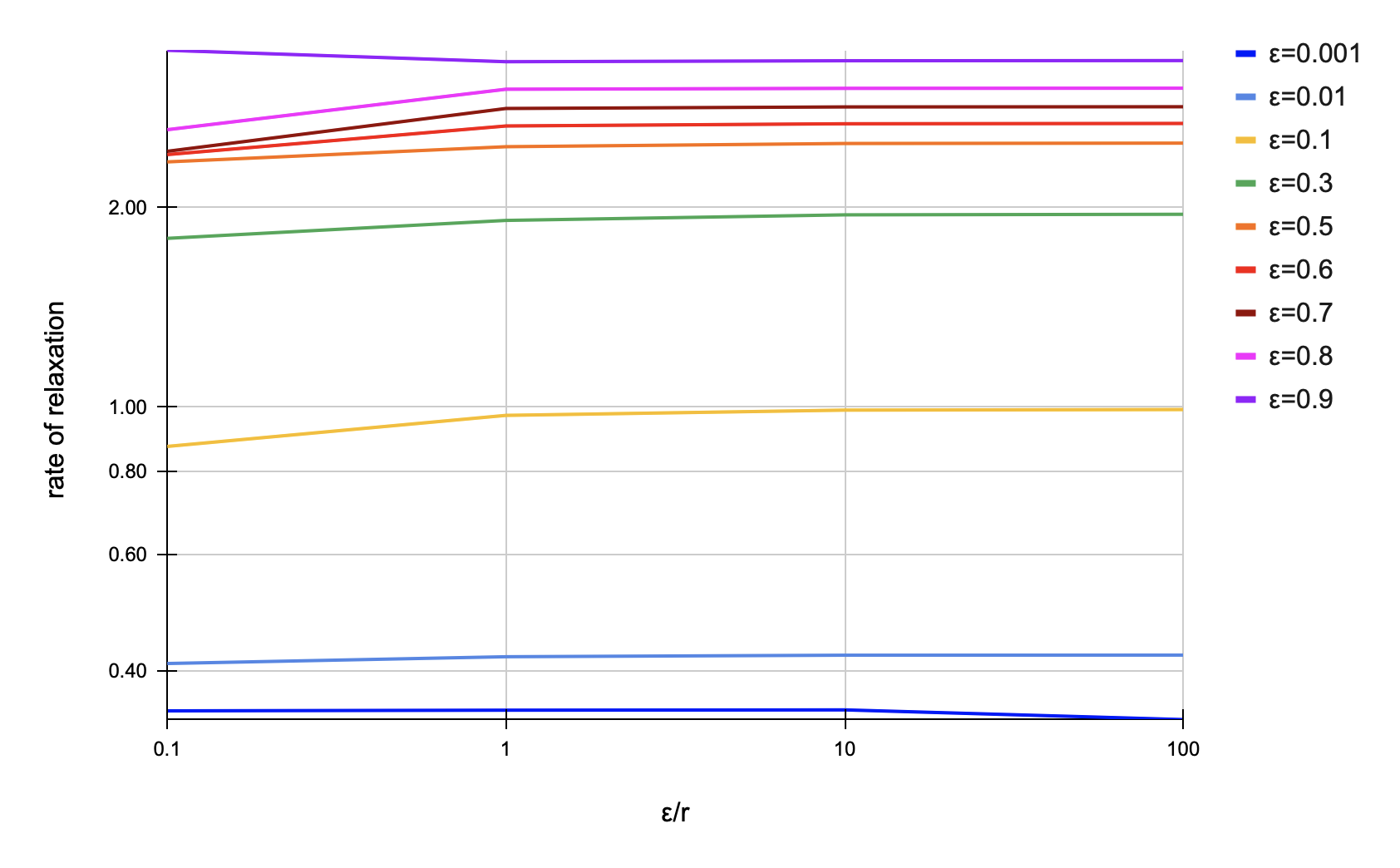}
    \caption{Relaxation time of the flattening of $u$} 
        \label{fig:flatten}
\end{figure}
\vincent{To distinguish between the diffusive and logistic regimes, we define three time scales: a ``diffusive time-scale," a ``transient time-scale", and a ``logistic time-scale". The diffusive time-scale, $t_D$, is defined as the first time for which the $L^2$-loss of $u_x$ relative to zero is small, that is, when $\|u_x(t_{D})\|^2< 0.01$. In this case, we say that $u$ has ``relaxed to the initial average" and refer to $t_D$ is the \textit{relaxation time for the initial average}. The corresponding relaxation rate is then defined simply as $t_D^{-1}$. Analogously, the ``transient time-scale" is defined as the first time, $\tau$, after $t_D$, that $u$ leaves an $L^2$-neighborhood of the initial average, that is, when  $\|u(\tau)-\intol u_0\|^2\geq 0.01$, where $\tau>t_D$. Lastly, the ``logistic time-scale" is the first time, $t_L$, after $\tau$, such that $u$ reaches a neighborhood of its carrying capacity, that is $\|u(t_L)-1\|^2<0.01$.}

\vincent{\cref{fig:ini:ave} shows that the relaxation time for the initial average seems to be largely independent of $r$, for each fixed $\veps$, indicating that the dominant dependence is on $\veps$.  However, when instead $r$ is fixed and $\veps$ is allowed to vary, \cref{fig:relax:ave} exposes a more nuanced dependence of the relaxation rate for the initial average on the ratio $\veps/r$. In particular, \cref{fig:relax:ave} shows that for each fixed $r$, the relaxation rate for the initial average depends increasingly on $\veps/r$. Moreover, as $r$ decreases to $0$, this dependence is maintained, indicating that the dominant dependence is indeed on the ratio $\veps/r$.
It is important to note that for $r\geq 1$ and small ratios of $\veps/r$, we observed that the solution $u$ will converge towards its carrying capacity in such a way that the ``diffusive regime" is bypassed, indicating a  bifurcation \textit{in the transient dynamics} for values of $r$ sufficiently large.} 

\vincent{An interesting and unexpected observation from our numerical experiments is the remarkably robust \textit{transition behavior} exhibited by the solutions across a wide range of parameter regimes. Indeed, \cref{fig:rel:time} and \cref{fig:table} show that the population density, as represented by $u$, tends to spend \vincent{\textit{a fixed percentage of time}, about $75\%$, between leaving the diffusive regime, that is, leaving a neighborhood of its initial average, and reaching the logistic regime, that is, reaching a neighborhood of its carrying capacity. This implies that the first $25\%$ of the time is spent near its initial average. The relative time can be computed from \cref{fig:table} by taking the difference between the relaxation time for the carrying capacity, $t_D$, and the transient time scale, $\tau$, then normalizing by $t_D$. We emphasize that in each of these experiments, the error thresholds defining the time-scales are \textit{maintained} across all choices of the parameters. While a rigorous proof of this may be difficult to obtain, the authors believe this transitory phenomenon is worth further investigation. In particular, though it was not investigated here, it would be interesting to carry out further numerical studies to test the robustness of this phenomenon across different initial data.} 
}

\vincent{Lastly, we also study the spatial homogenizing property of $u$ that was robustly observed in our numerical experiments. In particular, it was observed that $u$ tends to become constant in space before it reaches its carrying capacity (see \cref{fig:transient:b} and \cref{fig:logistic}). 
Since the diffusivity of the bacterial population is exactly $1$, this means that the chemical diffusivity is always smaller than the diffusivity of the bacterial population. It is, thus, interesting to see if there is any change in the behavior of the solutions when the coefficients are of the same order. In particular, we study this behavior for values of $\veps$ approaching $1$. We track this by viewing $t_D^{-1}$ as a ``rate of flattening" of $u$. Then \cref{fig:flatten} shows that when $\veps/r\geq 1$, flattening rate becomes constant, and is increasing as $\veps$ increases. The dependence in $r$ is largely absent when the ratio $\veps/r\geq1$, but expectedly comes into effect when $\veps/r<1$. An elementary heuristic from rescaling seems to explain this phenomenon adequately. Indeed, by rescaling time according to  $t'=\frac tr$, one obtains that the diffusive coefficient of $u$ becomes $\frac{1}{r}$ and the one of $v$ becomes $\frac{\veps}{r}$. Thus, at least in logistic time scales, the situation when $\veps\sim r$ and $r\gtrsim1$ can be viewed as a chemical diffusion-dominant regime, so that the rate of relaxation to a constant, that is the tendency of $u$ to spatially homogenize via diffusive effects becomes slower as $\veps/r\rightarrow0$.}


\section{Conclusion and future work}
\vincent{We have shown that the logistic growth for this particular model for chemotaxis over the unit interval dominates the long time behavior of the solutions in the case that the population density $u$ satisfies zero flux boundary conditions and the chemical concentration, insofar as it is represented through $v$, satisfies Dirichlet boundary conditions. In particular, all solutions emanating from a large class of initial data in $H^2(0,1)$ respecting the boundary conditions, necessarily have the population density converging to the carrying capacity, while the chemical density approaches a constant. We also show that the solutions of the diffusive system converge in the $H^1$--norm towards solutions with zero-diffusivity, therefore, verifying that the system with chemical diffusion can serve as suitable approximation to the non-diffusive system. Lastly, numerical investigations are carried out that verify the rigorous results and, additionally, shed light on the transient behavior of the solutions. With regard to the latter, several numerical experiments indicate interesting, non-trivial properties of the transient dynamics that warrant further mathematical and numerical investigation.} 
\section{Acknowledgements}
The authors would like to thank Dr. Ricardo Cortez for the guidance in the numerical simulations, as well as for the insightful discussions. The work
of Vincent R. Martinez was partially supported by the PSC-CUNY
Research Award Program under grant PSC-CUNY 62239-00 50. The research of K. Zhao was partially supported by the Simons Foundation Collaboration Grant for Mathematicians No. 413028.

\bibliographystyle{alpha}

\end{document}